\documentclass{daj}

\dajAUTHORdetails{%
  title = {Approximate Invariance for Ergodic Actions of Amenable Groups}, 
  author = {Michael Bj\"orklund, Alexander Fish},
  plaintextauthor = {Michael Bjoerklund, Alexander Fish},
  plaintexttitle = {Approximate Invariance for Ergodic Actions of Amenable Groups}, 
  runningtitle = {Approximate Invariance for Ergodic Actions}, 
  runningauthor = {Michael Bj\"orklund, Alexander Fish},
  copyrightauthor = {Michael Bj\"orklund and Alexander Fish},
  keywords = {Action sets, aperiodicity, density theorems},
}   

\dajEDITORdetails{%
   year={2019},
   number={6},
   received={3 October 2018},   
   published={24 May 2019},  
   doi={10.19086/da.8471},       
}   

\usepackage[latin1]{inputenc}
\usepackage{tikz}
\usetikzlibrary{shapes,arrows}
\usetikzlibrary{arrows.meta}
\usepackage{forest}
\usepackage{tikz-qtree}

\usetikzlibrary{arrows,shapes,positioning,shadows,trees}

\usepackage{etoolbox}
\usepackage{amsmath}
\usepackage{enumerate}
\usepackage{amssymb}
\usepackage{amscd}
\usepackage{amsthm}
\usepackage{amsfonts}
\usepackage{graphicx}
\usepackage[all,cmtip]{xy}
\usepackage{enumitem}

\patchcmd{\subsection}{-.5em}{.5em}{}{}
\patchcmd{\subsubsection}{-.5em}{.5em}{}{}

\usepackage{enumitem}

\usepackage[T1]{fontenc}

\usepackage[scaled]{helvet} 
\usepackage{courier} 
\usepackage{eulervm}
\normalfont

\makeatother
\usepackage{hyperref}

\numberwithin{equation}{section}

\newcommand{\im}{\operatorname{Im}}

\newcommand{\cB}{\mathcal{B}}

\newcommand{\cE}{\mathcal{E}}
\newcommand{\cF}{\mathcal{F}}

\newcommand{\cJ}{\mathcal{J}}
\newcommand{\cK}{\mathcal{K}}
\newcommand{\cL}{\mathcal{L}}
\newcommand{\cM}{\mathcal{M}}

\newcommand{\cP}{\mathcal{P}}

\newcommand{\bC}{\mathbb{C}}

\newcommand{\bE}{\mathbb{E}}

\newcommand{\bN}{\mathbb{N}}

\newcommand{\bQ}{\mathbb{Q}}
\newcommand{\bR}{\mathbb{R}}

\newcommand{\bT}{\mathbb{T}}

\newcommand{\bZ}{\mathbb{Z}}

\newcommand{\ra}{\rightarrow}

\newcommand{\qor}{\quad \textrm{or} \quad}

\newcommand{\qand}{\quad \textrm{and} \quad}

\def\acts{\curvearrowright}
\newcommand\subsetsim{\mathrel{%
\ooalign{\raise0.2ex\hbox{$\subset$}\cr\hidewidth\raise-0.8ex\hbox{\scalebox{0.9}{$\sim$}}\hidewidth\cr}}}
\newcommand{\eps}{\varepsilon}

\DeclareMathOperator{\supp}{supp}

\theoremstyle{theorem}
\newtheorem{theorem}{Theorem}[section]
\newtheorem{corollary}[theorem]{Corollary}
\newtheorem{proposition}[theorem]{Proposition}
\newtheorem{lemma}[theorem]{Lemma}

\theoremstyle{definition}
\newtheorem{definition}[theorem]{Definition}
\newtheorem{remark}[theorem]{Remark}

\tikzstyle{decision} = [diamond, draw, fill=blue!20, 
    text width=4.5em, text badly centered, node distance=3cm, inner sep=0pt]
\tikzstyle{block} = [rectangle, draw, fill=blue!20, 
    text width=5em, text centered, rounded corners, minimum height=2em]
\tikzstyle{line} = [draw, -latex']
\tikzstyle{cloud} = [draw, ellipse,fill=red!20, node distance=3cm,
    minimum height=2em]

\renewcommand\labelenumi{(\roman{enumi})}
\renewcommand\theenumi\labelenumi

\usepackage{changepage}

\begin{document}

\begin{frontmatter}[classification=text]

\title{Approximate Invariance for Ergodic Actions of Amenable Groups}

\author[michael]{Michael Bj\"orklund\thanks{Supported by European Framework Program (FP7/2007/2012 Grant Agreement 203418) when
he was a postdoctoral fellow at Hebrew University, ETH Fellowship FEL-171-03 between January 2011 and August 2013, and 
GoCAS (Gothenburg Centre for Advanced Studies) since September 2013.}}
\author[sasha]{Alexander Fish}

\begin{abstract}
We develop in this paper some general techniques to analyze action sets of 
small doubling for probability measure-preserving actions of amenable groups.

As an application of these techniques, we prove a dynamical generalization of 
Kneser's celebrated density theorem for subsets in $(\bZ,+)$, valid for any countable
amenable group, and we show how it can be used to establish a plethora of new inverse 
product set theorems for upper and lower asymptotic densities. We provide several examples 
demonstrating that our results are optimal for the settings under study. 
\end{abstract}
\end{frontmatter}


\tableofcontents
\setcounter{tocdepth}{2}

\addtocontents{toc}{\protect\setcounter{tocdepth}{1}}
\section{Introduction}

\subsection{Motivation}

Subsets of locally compact groups which are almost closed under multiplication are classical 
objects of study in harmonic analysis, number theory and geometry, and continually appear in new applications. 
Despite the fact that these objects are often far from being actual subgroups, they nevertheless seem to obey  some 
form of approximate group theory, whose foundation is still in its early infancy. With this paper we wish to 
take the first steps towards what could be called approximate dynamics or approximate ergodic theory, 
a line of research concerned with the interplay between expansion and approximate invariance of so called action sets, 
dynamical analogues of product sets in groups. 

\begin{definition}
Let $G$ be a group and let $Y$ be a set upon which $G$ acts. Given $A \subset G$ and $B \subset Y$, their
\emph{action set} $AB$ is defined by
\[
AB = \bigcup_{a \in A} aB = \big\{ ab \, : \, a \in A, \enskip b \in B \big\} \subset Y.
\]
Let $k \geq 1$ be an integer. We say that $(A,B)$ is \emph{$k$-doubling} if there exists a finite subset $F \subset G$ 
of cardinality at most $k$ such that $AB \subset FB$; in this case, we also say that $B$ is 
\emph{$k$-approximately invariant} under $A$. To avoid trivialities, one usually insists that $k$ is "small" compared to
the "sizes" of the sets $A$ and $B$.

In the case when $Y = G$, endowed with the natural $G$-action on itself from the left, and 
$A,B \subset G$, then the action set $AB$ is called the \emph{product set} of $A$ and $B$, and $A \subset G$ is 
called \emph{$k$-doubling} if $(A,A)$ is a $k$-doubling pair. We say that $A \subset G$ is a 
\emph{$k$-approximate subgroup} if it is $k$-doubling, symmetric and contains the identity element $e_G$. 
\end{definition}

Let us begin by making a trivial observation: If $G \acts Y$, then $k$-doubling sets in $G$ naturally give 
rise to $k$-approximately invariant subsets in $Y$ as follows. Suppose that $A \subset G$ is $k$-doubling and $B_o \subset Y$ 
is any subset; then $B = A B_o$ is $k$-approximately invariant under $A$ - indeed, since $A^2 \subset FA$ for some 
subset $F \subset G$ with $|F| \leq k$, we have
\[
AB = A^2 B_o \subset FA B_o = FB.
\]
A fundamental line of research is to investigate to which extent the converse holds, that is to say, if a set $B \subset Y$
is $k$-approximately invariant under a subset $A \subset G$, must then $A$ be contained in a $k'$-approximate 
subgroup of $G$ where $k'$ is not much larger than $k$. If $k = 1$, i.e. if $B$ is invariant (on the nose) under $A$, then $B$
is clearly also invariant under the subgroup of $G$ which is generated by $A$. The case $k = 2$ is already much harder to deal
with, and will be the focus of our investigations here. \\

Before we outline the theme of this paper, and state our main results, we say a few words about history. The term "$k$-approximate subgroup" was coined by Tao in \cite{ta08}, but implicit uses of the notion can be 
traced back much further. For instance, the study of $2$-approximate subgroups was initiated in the work by 
Mann \cite{ma1942} on Schnirelmann densities of sumsets in $(\bZ,+)$, which Khintchine \cite{ch1951} later 
referred to as one of the "Three Pearls of Number Theory", and it was continued in the subsequent work by Kneser 
\cite{kn1953}, as well as in many important works by Kneser \cite{kn1956}, Kempermann \cite{ke1964} and 
others on product sets in compact groups. 

However, the impetus to Tao's work was the early works of Freiman \cite{fr73} on general \emph{finite} $k$-approximate 
subgroups in $(\bZ,+)$, later extended by Ruzsa, as well as in the more recent works by Bourgain-Gamburd \cite{boga2008} 
and Helfgott \cite{he08} on finite $k$-approximate subgroups in finite simple groups. A very general theorem in this direction 
was recently established by Breuillard, Green and Tao \cite{brgrta2012}.

From a very different point of view, Yves Meyer \cite{me1972} began in the sixties his very influential study of "large" and discrete 
$k$-approximate subgroups in Euclidean spaces, which are today mostly known under the names \emph{quasicrystals}, 
\emph{Meyer sets} or \emph{approximate lattices}. Extensions to non-abelian groups were recently developed by the first author
and Tobias Hartnick \cite{bjha17, bjha172,bjha173}. \\

As a warm-up, we provide a classification of $2$-doubling pairs $(A,B)$, where 
$A$ is a "large" and "aperiodic" (or "spread-out") subset of a countable (infinite) abelian group $G$, for instance $(\bZ,+)$, and $B$ is a Borel set with 
positive measure in some ergodic Borel $G$-space  $(Y,\nu)$. In other words, we shall assume that there are $s,t \in G$ such 
that
\[
AB \subset sB \cup tB.
\]
To avoid trivialities, we wish to exclude the case when $sB \cup tB = Y$, which can be done by assuming that $\nu(B) < 1/2$.
If one further assumes that $A$ is "bigger" than $B$, then we wish to show that $A$ must "essentially" coincide with a $2$-approximate subgroup 
of a very special form. 
\begin{definition}[Large/Spread-out]
Let $G$ be a countable amenable group and let $d^*$ denote the upper Banach density on $G$ (see \eqref{defBanachdens} for
the definition). We say that $A \subset G$ is 
\begin{itemize}
\item \emph{large} if $d^*(A) > 0$.
\item \emph{spread-out} if $A$ is large and $G_o A_o = G$ for every finite index subgroup $G_o$ in $G$ and for every $A_o \subset A$
having $d^*(A_o) = d^*(A)$. 
\end{itemize}
In other words, $A$ is spread-out if one cannot pass to a subset with the same upper Banach density which is contained in a proper periodic subset of $G$. In particular, $A$, as well as any of its subsets with the same upper 
Banach density, projects onto every finite quotient of $G$. We note that if $G$ lacks proper finite-index subgroups, for instance if $G = (\bQ,+)$, then every large set is automatically spread-out.
\end{definition}
\begin{remark}
Our definition of a spread-out set might be somewhat hard to digest, and perhaps it seems to be a bit too strong of an assumption; 
we have chosen the formulation above to make certain parts of our arguments run smoother, but it will be clear from our proofs that 
one needs much less. For instance, one objection to the phrasing above could be that the condition $d^*(A_o) = d^*(A)$ is not very informative; with some additional work, one could prove that all of our main results still hold if one in addition insists that this identity is realized along the \emph{same} F\o lner sequence (see Subsection \ref{subsection:amenable} for definitions). However, to keep the exposition clean,
we shall refrain from such technical indulgences. 
\end{remark}

\begin{definition}[Group compactifications and induced actions]
Let $G$ be a countable group and let $K$ be a compact and second countable group with Haar probability measures $m_K$. Suppose that there 
exists a homomorphism $\tau : G \ra K$ with dense image. We can then endow $K$ with a continuous $m_K$-preserving $G$-action by 
\begin{equation}
\label{defKtau}
g.k = \tau(g)k \quad \textrm{for $g \in G$ and $k \in K$}.
\end{equation}
With this notation understood, we give the following definitions.
\begin{itemize}
\item The pair $(K,\tau)$ is called a \emph{group compactification} of $G$.
\item The Borel $G$-space $(K,m_K)$, with the $G$-action defined above, is called the \emph{induced $G$-space} associated to the group
compactification $(K,\tau)$.
\item If $(Y,\nu)$ is a Borel $G$-space and $Y_o \subset Y$ is a $G$-invariant $\nu$-conull subset, then a $G$-equivariant Borel map $\sigma : Y_o \ra K$,
where $K$ is endowed with the $G$-action above, is called a \emph{$G$-factor map}. The dependence on the $\nu$-conull subset $Y_o$ will be suppressed and we shall denote the $G$-factor map by $\sigma : (Y,\nu) \ra (K,m_K)$.
\end{itemize}
\end{definition}

The exact formulation of our classification now reads as follows.  

\begin{theorem}[Warm-up]
\label{warm-up}
Let $G$ be a countable abelian group and let $G \acts (Y,\nu)$ be a totally ergodic Borel $G$-space. Suppose that 
$A \subset G$ is spread-out, and $B \subset Y$ is a Borel set with positive $\nu$-measure such that
\begin{equation}
\label{assThmA}
\nu(B) \leq d^*(A) < 1/2 \qand AB \subset sB \cup tB, \quad \textrm{modulo $\nu$-null sets}, 
\end{equation}
for some $s,t \in G$. Then $d^*(A) = \nu(B)$, and there exist
\begin{enumerate}
\item a torus compactification $(\bT,\tau)$ of $G$ and a closed interval $I_o \subset \bT$ with $m_{\bT}(I_o) = d^*(A)$,
\item a $G$-factor map $\sigma : (Y,\nu) \ra (\bT,m_{\bT})$ and a closed interval $J_o \subset \bT$ with $m_{\bT}(J_o) = \nu(B)$,
where the $G$-action on $\bT$ is defined as in \eqref{defKtau} using the group compactification $(\bT,\tau)$,
\end{enumerate}
such that $A \subset \tau^{-1}(I_o)$ and $B = \sigma^{-1}(J_o)$ modulo $\nu$-null sets.
\end{theorem}

The theorem shows that the structure of $2$-doubling pairs for ergodic actions is very rigid; the set $B$ must stem from an interval
in one-dimensional torus, and the set $A$ is contained in a set $S$ of the form $\tau^{-1}(I_o)$, where $I_o$ is a closed interval in
the same one-dimensional torus. Furthermore, $S$ has the same upper Banach density as $A$ (this follows for instance from 
Corollary \ref{cor_bohrset} in the appendix). It is straightforward to check that $S$ is a $3$-approximate subgroup,
and in fact a $2$-approximate subgroup if the endpoints of the interval $I_o$ belong to $\tau(G)$. We stress that the converse also
holds; if $A$ and $B$ are as in the conclusion of Theorem \ref{warm-up}, then $(A,B)$ is $2$-doubling (modulo null sets).

\subsection{Main dynamical results}

Let us now connect Theorem \ref{warm-up} to the main theme of this paper. If $A$ and $B$ are as in this theorem, then it is clear
that
\[
\nu(AB) \leq 2 \nu(B) \leq d^*(A) + \nu(B) < 1.
\]
Theorem \ref{erg1} below tells us that if $A$ is a spread-out subset, then the reverse inequality $\nu(AB) \geq \min(1,d^*(A) + \nu(B))$ 
holds for \emph{all} Borel sets $B \subset Y$, so in the setting at hand, we must have
\[
\nu(B) = d^*(A) \qand \nu(AB) = d^*(A) + \nu(B) < 1.
\]
Theorem \ref{warm-up} now follows from the latter part of Theorem \ref{erg2} below. \\

The general framework for the theorems below reads as follows:
\begin{itemize}
\item $G$ - a countable amenable group.
\item $d^*$ - the upper Banach density on $G$.
\item $(Y,\nu)$ - a standard Borel probability measure space, equipped with an ergodic action of $G$ by measure-preserving
bijections.
\end{itemize}

Our first main result (which is proved in Subsection \ref{subsec:maindyn}) asserts that if an action set $AB$ in $Y$ is "small"
with respect to $d^*(A)$ and $\nu(B)$, then $A$ must be "close" to a periodic subset.

\begin{theorem}
\label{erg1}
Let $A \subset G$ be a large set and $B \subset Y$ a Borel set with positive measure.
Suppose that either
\begin{enumerate}
\item $A$ is spread-out \textbf{or}
\item all finite quotients of $G$ are \textsc{abelian}, and there is no finite-index subgroup $G_o < G$
such that
\begin{equation}
\label{perbnd}
G_o A \neq G \qand d^*(G_o A) < d^*(A) + \frac{1}{[G : G_o]},
\end{equation}
\end{enumerate}
then $\nu(AB) \geq \min\big(1,d^*(A) + \nu(B)\big)$.
\end{theorem}

\begin{remark}
If $G$ is abelian, then every finite quotient group is of course also abelian. For a non-abelian example of $G$ for which every
finite quotient group is abelian, consider
the lampligher group $G = \bQ \wr \bZ$ (with $\bQ$-valued "lamps"). This is the solvable (hence amenable) wreath 
product of the two abelian 
groups $(\bQ,+)$ and $(\bZ,+)$, and it is not hard to show that every homomorphism of $G$ onto a 
finite group factors through $\bZ$, and thus every finite factor group of $G$ is abelian. 
\end{remark}

Let us now try to understand when the lower bound on $\nu(AB)$ in the previous theorem is attained. We saw in
Theorem \ref{warm-up} that if $G$ is abelian, then a special role is played by sets of the form $\tau^{-1}(I_o)$, 
where $(\bT,\tau)$ is a torus compactification of $G$ and $I_o$ a closed interval in $\bT$. These are examples of 
so called \emph{Sturmian sets}, which we now define for general groups.

\begin{definition}[Sturmian set]
\label{def:sturm}
Let $G$ be a countable group and let $M$ denote either $\bT$ or $\bT \rtimes \{-1,1\}$, where in the latter case, 
$\{-1,1\}$ acts by multiplication on $\bT$. We say that a subset $S \subset G$ is \emph{Sturmian} if there exist
\begin{itemize}
\item a homomorphism $\tau : G \ra M$ with dense image, 
\item a closed symmetric interval $I_o \subset \bT$ and $t_o \in M$,
\end{itemize}
such that either $S = \tau^{-1}(I_o t_o)$ if $M = \bT$ (in which case the assumption that $I_o$ is symmetric can be dropped) 
or $S = \tau^{-1}((I_o \rtimes \{-1,1\})t_o)$ if $M = \bT \rtimes \{-1,1\}$. In the latter case, we say that $S$ is \emph{twisted}. 
\end{definition}

Clearly, abelian groups do not admit twisted Sturmian sets. On the other hand, one can readily check that the infinite
dihedral group \emph{only} admits twisted Sturmian sets, and no "un-twisted" ones. Sturmian sets in $\bZ^d$ have 
been extensively studied in complexity theory and tiling theory, see for instance the survey \cite{wi2007}, while we seem to
be the first to address their twisted analogues. \\

If $(M,\tau)$ is as in Definition \ref{def:sturm}, then we get an ergodic Borel $G$-space by letting $G$ act on $M$ by 
multiplication on the left via $\tau$. This action clearly preserves the Haar measure $m_M$ on $M$ and is ergodic (because the
image of $\tau$ is dense, see for instance Lemma \ref{basic_isomG} below). Let $I_o$ and $J_o$ be symmetric closed 
intervals of $\bT$ with $m_{\bT}(I_o + J_o) < 1$, and set
\[
A = \tau^{-1}(I_o) \subset G \qand B = J_o \subset M,
\]
if $M = \bT$ (or the twisted versions if $M = \bT \rtimes \{-1,1\}$). Then it is not hard to show that 
\begin{itemize}
\item $A$ is spread out and $d^*(A) = m_{\bT}(I_o)$,
\item $AB$ does not contain, modulo $m_M$-null sets, a Borel set with positive measure
which is invariant under a finite-index subgroup of $G$, and
\item $m_M(AB) = d^*(A) + m_M(B) < 1$.
\end{itemize}

Our second main theorem (which is proved in Subsection \ref{subsec:maindyn}) asserts that upon passing to factors, the setting described above is the only source of examples 
of ergodic Borel $G$-spaces for which the lower bound in Theorem \ref{erg1} is attained (this is spelled out precisely for 
abelian groups below, but the version for general amenable groups will be clear from the proofs).

\begin{theorem}
\label{erg2}
Let $A \subset G$ be a large set and $B \subset Y$ a Borel set with positive measure.
Suppose that
\begin{enumerate}
\item $A$ is spread-out,
\item $AB$ does not contain, modulo $\nu$-null sets, a Borel set with positive measure
which is invariant under a finite-index subgroup of $G$,
\end{enumerate}
and
\[
\nu(AB) = d^*(A) + \nu(B) < 1.
\]
Then $A$ is contained in a Sturmian set with the same upper Banach density as $A$. \\

If $G$ is abelian and $G \acts (Y,\nu)$ is totally ergodic, then Condition (ii) is automatic, and there exist
\begin{enumerate}
\item a torus compactification $(\bT,\tau)$ and closed intervals $I_o, J_o \subset \bT$ with 
\[
d^*(A) = m_{\bT}(I_o) \qand A \subset \tau^{-1}(I_o) \qand \nu(B) = m_{\bT}(J_o).
\]
\item a $G$-factor map $\sigma : (Y,\nu) \ra (\bT,m_{\bT})$ such that $B = \sigma^{-1}(J_o)$ modulo $\nu$-null sets, where 
$G$ acts on $\bT$ via $\tau$ as in \eqref{defKtau}.
\end{enumerate}
\end{theorem}

\begin{remark}
\label{perfect1}
The first part of Theorem \ref{erg2} has a curious consequence: if there are sets $A$ and $B$ as above, then $G$ must admit a 
non-trivial homomorphism into either $\bT$ or $\bT \rtimes \{-1,1\}$. In particular, if $G$ is perfect, that is to say, if $G = [G,G]$,
then we cannot find a spread-out set $A \subset G$ and $B \subset Y$ as above which satisfy the identity $\nu(AB) = d^*(A) + \nu(B) < 1$.

It is natural to ask how "close" to this case we can get if $G$ is perfect. Corollary \ref{addcons2} in Subsection \ref{subsec:aux}
shows that for some perfect amenable groups (for instance, the Grigorchuk group), the condition that $A$ is spread-out 
is so strong that it forces $\nu(AB) = 1$ for every Borel set $B \subset Y$ with positive $\nu$-measure. Similarly, if $G$ is an
amenable weakly mixing (minimally almost periodic) group, then $\nu(AB) = 1$ whenever $A$ is large and $B \subset Y$ has 
positive $\nu$-measure. Since every (non-trivial) quotient of a weakly mixing group is weakly mixing, we readily see that weakly mixing groups 
must be perfect (as every abelian group admits actions which are not weakly mixing). 
\end{remark}

\subsection{Applications to product sets}

In this section we translate our dynamical results into density combinatorial results. While this is rather straightforward to 
do for \emph{Banach densities}, we need to develop new tools in order to address \emph{asymptotic densities}. This is due
to the fact that when dealing with asymptotic densities from a "dynamical" point of view, one needs to use action set theorems 
for \emph{non-ergodic} Borel $G$-spaces as well, and our dynamical results above only apply to \emph{ergodic} ones. 

Starting from Subsection \ref{subsec:dealwithit}, we describe how one can transfer inverse theorems for action set with respect
to \emph{ergodic} Borel $G$-spaces to inverse theorems for action sets with respect to \emph{non-ergodic} Borel $G$-spaces. On
the way, we shall need some terminology which we now introduce.

\begin{definition}[Thick/piecewise periodic/syndetic]
We say that a subset $P \subset G$ is
\begin{itemize}
\item \emph{thick} if for every finite set $F \subset G$, there exists $g \in G$ such that $Fg \subset P$.
\item \emph{periodic} if there exists a finite-index subgroup $G_o < G$ such that $G_o P = P$.
\item \emph{piecewise periodic} if there exist a periodic set $Q$ and a thick set $T$ such that $P = Q \cap T$.
\item \emph{syndetic} if it intersects non-trivially every thick set in $G$, or equivalently, if there exists a finite set $F$ in $G$
such that $FP = G$.
\end{itemize}
\end{definition}

\subsubsection{Inverse theorems for Banach densities}

Let us begin our discussions here with a classification of spread-out $2$-approximate subgroups of countable 
amenable groups. 

\begin{theorem}[Warm-up]
\label{warmup2}
Suppose that $A \subset G$ is a spread-out $2$-approximate subgroup such that $A^2$ does
not contain a piecewise periodic set. Then $A$  is contained in a Sturmian subset of $G$ with 
the same upper Banach density as $A$. 

If $G$ does not have any proper finite-index subgroups, then it suffices to assume that $A$ is 
large and $A^2$ is not thick.
\end{theorem}

\begin{remark}
We observe, just as we did in Remark \ref{perfect1}, that this theorem in particular implies that 
amenable perfect groups do not admit \emph{spread-out} $2$-approximate subgroups (whose
squares do not contain piecewise periodic sets). In particular, perfect amenable groups without 
proper finite index subgroups do not admit spread-out $2$-approximate subgroups with upper Banach
densities strictly between $0$ and $1/2$. 
\end{remark}

Theorem \ref{warmup2} is a rather direct consequence of Theorem \ref{dens1ubd} and Theorem 
\ref{dens2ubd} below. Indeed, if $A$ is a $2$-approximate subgroup with $d^*(A^2) < 1$, then 
\[
d^*(A^2) \leq 2d^*(A).
\]
If $A$ is spread-out, then Theorem \ref{dens1ubd} will show that the reverse inequality also holds,
and thus $d^*(A^2) = 2 d^*(A)$. If $A^2$ does not not contain a piecewise periodic set, then Theorem
\ref{dens2ubd} now implies that $A$ is contained in a Sturmian set with the same upper Banach density,
finishing the proof of Theorem \ref{warmup2}. \\

The arguments needed to prove the following two theorems from Theorem \ref{erg1} and Theorem 
\ref{erg2} are nowadays standard, and are only sketched in the beginning of Subsection \ref{subsec:maindens}.

\begin{theorem}
\label{dens1ubd}
If $A \subset G$ is spread-out, then 
\[
d^*(AB) \geq \min(1,d^*(A) + d^*(B)), \quad \textrm{for all large $B \subset G$},
\]
and
\[
d_*(AB) \geq \min(1,d^*(A) + d_*(B)), \quad \textrm{for all syndetic $B \subset G$}.
\]
If one assumes that every finite quotient of $G$ is abelian, then, instead of assuming
that $A$ is spread-out, 
it suffices to assume that there is no finite-index subgroup $G_o < G$ such that 
\eqref{perbnd} holds for $A$.
\end{theorem}

\begin{remark}
A version of this theorem for $G = (\bN,+)$ was proved by Jin \cite{ji2010} using very
different methods. Griesmer \cite{gr2013} proved a version of Theorem \ref{dens1ubd}
(as well as of Theorem \ref{dens2ubd} below) for countable abelian groups; his proof was very much inspired by some earlier versions
of our Correspondence Principles for product sets (see Proposition \ref{corr1} and 
Proposition \ref{corr2} below).
\end{remark}

\begin{theorem}
\label{dens2ubd}
Let $A \subset G$ be spread-out and $B \subset G$ large and suppose that $AB$ does not contain a 
piecewise periodic set. If $d^*(AB) = d^*(A) + d^*(B) < 1$, then $A$ is contained in a Sturmian set with
the same upper Banach density as $A$.
\end{theorem}

\subsubsection{Inverse theorems for asymptotic densities}

Our inspiration for this subsection comes from a classical result of Kneser \cite{kn1953}, generalizing an earlier 
landmark made by Mann \cite{ma1942}. It is an inverse result for subsets $A, B \subset \bN$ with positive lower asymptotic densities along
the F\o lner sequence $([1,n])$ such that 
\[
\underline{d}_{[1,n]}(A+B) < \min(1,\underline{d}_{[1,n]}(A) + \underline{d}_{[1,n]}(B)),
\]
and roughly asserts that $A$ and $B$ must be contained in periodic subsets which are not "much larger" in size than
$A$ and $B$, and moreover, $A+B$ is a co-finite subset of a periodic set in $\bN$. In particular, if $A$ is sufficiently
"aperiodic", then 
\[
\underline{d}_{[1,n]}(A+B) \geq \min(1,\underline{d}_{[1,n]}(A) + \underline{d}_{[1,n]}(B))
\] 
for all $B \subset \bN$ with positive lower asymptotic density along $([1,n])$. 

Our aim here is to generalize the latter (weaker) formulation of Kneser's Theorem to a general countable amenable 
group $G$ and to a general F\o lner sequence $(F_n)$ therein. We shall also prove an inverse theorem in the case 
when the lower bound is attained. We stress that our results are new already in the case $G = (\bZ,+)$ and 
$F_n = [1,n]$, and in Appendix \ref{sec:pec} we provide several examples showing that our results are optimal in the
case $G = (\bZ,+)$ and $F_n = [-n,n]$. \\

In what follows, let $G$ be a countable amenable group. Our first result improves (albeit under stronger conditions on
$A$ and $B$), the lower bound in the (weak) formulation of Kneser's Theorem above, by replacing $\underline{d}_{([1,n]}(A)$
with the (possibly) larger quantity $d^*(A)$. The proof can be found in Subsection \ref{prfdens1}.

\begin{theorem}
\label{dens1}
Let $(F_n)$ be a F\o lner sequence in $G$. Suppose that 
\begin{enumerate}
\item $A \subset G$ is spread-out.
\item $B \subset G$ is syndetic.
\item $AB$ is not thick.
\end{enumerate}
Then, 
\begin{equation}
\label{coolbnds}
\overline{d}_{(F_n)}(AB) \geq d^*(A) + \overline{d}_{(F_n)}(B)
\qand
\underline{d}_{(F_n)}(AB) \geq d^*(A) + \underline{d}_{(F_n)}(B).
\end{equation}
If every finite quotient of $G$ is abelian, then instead of (i), we only need to assume that there is no finite-index 
subgroup $G_o < G$ such that \eqref{perbnd} holds.
\end{theorem}

\begin{remark}
In Proposition \ref{attempt1} and Proposition \ref{attempt2} we show that Condition (iii) and Condition (ii) respectively 
cannot be dispensed with, already in the case $G = (\bZ,+)$ and $F_n = [-n,n]$. 
\end{remark}

Our second theorem addresses the equality case in Theorem \ref{dens1}. 

\begin{theorem}
\label{dens2}
Let $(F_n)$ be a F\o lner sequence in $G$. Suppose that 
\begin{enumerate}
\item $A \subset G$ is spread-out.
\item $B \subset G$ is syndetic.
\item $AB$ does not contain a piecewise periodic subset.
\item Either 
\[
\overline{d}_{(F_n)}(AB) = d^*(A) + \overline{d}_{(F_n)}(B) < 1
\qor
\underline{d}_{(F_n)}(AB) = d^*(A) + \underline{d}_{(F_n)}(B) < 1.
\]
\end{enumerate}
Then $A$ is contained in a Sturmian set with the same upper Banach density as $A$.
\end{theorem}

\begin{remark}
In Proposition \ref{attempt3} and Proposition \ref{attempt4} we show that Condition (ii) and Condition (iii) respectively 
cannot be dispensed with, already in the case $G = (\bZ,+)$ and $F_n = [-n,n]$. 
\end{remark}

\subsection{Counterexamples}

Some readers might interpret the asymmetry in the roles of the sets $A$ and $B$ in Theorem \ref{dens1ubd} 
and Theorem \ref{dens2ubd} as a sign of incompleteness; after all, in these theorems, we only make assertions about $A$, and
say nothing about $B$. Of course, if $G$ is abelian, one can
simply swap the order of $A$ and $B$, and then use the theorems above to deduce things about the set $B$. However, there
is no reason why this trick should work if $G$ is non-abelian. 

The aim of the next two results (which are proved in Section \ref{sec:cntex}) is to show that if 
$G$ is sufficiently non-abelian (certain two-step solvable groups will do), then the roles of $A$ and $B$ are truly
asymmetric, and the conclusions about the set $A$ in Theorem \ref{dens1ubd} and Theorem \ref{dens2ubd} do \emph{not} hold for
$B$. Both results are derived from a general counterexample machine for semi-direct products which should be of independent interest.

\begin{theorem}
\label{thm_counterex1}
There is a countable two-step solvable group $G$ and $A \subset G$ with $d^*(A) = 1/2$
such that for every $0 < \eps < 1/2$, there is $B \subset G$ with $d^*(B) = \eps$, 
with the property that 
\[
d^*(A B) = d^*(A) < d^*(A) +d^*(B) < 1,
\]
and for every finite-index subgroup $G_o < G$, we either have
\[
G_o B = G \qor d^*(G_o B) > d^*(B) + \frac{1}{[G : G_o]}.
\]
\end{theorem}

\begin{theorem}
\label{thm_counterex2}
There is a countable two-step solvable group $G$ and $A, B \subset G$ such that
\begin{enumerate}
\item $A$ and $B$ are spread-out,
\item $A B$ does not contain a piecewise periodic set,
\end{enumerate}
and
\[
d^*(AB) = d^*(A) + d^*(B) < 1,
\]
but $B$ is \textsc{not} contained in a Sturmian set with the same upper Banach density.
\end{theorem}

\subsection{A few words about the proofs of the dynamical results}

To prove Theorem \ref{erg1} and Theorem \ref{erg2} we argue along the following lines. Let $G$ be a countable amenable
group and suppose that $G$ acts by homeomorphisms on a compact metrizable space $X$. Given a subset $A \subset X$
and $x \in X$, we write $A_x = \big\{ g \in G \, : \, gx \in A \big\}$, which provides us with a subset of $G$. It is not hard to prove
that \emph{every} subset of $G$ can be written in this form for \emph{some} compact $G$-space $X$, \emph{clopen} subset 
$A \subset X$ and base point $x_o \in X$ (see Subsection \ref{subsec:hull}). Let us fix such a triple $(X,x_o,A)$ and an ergodic
Borel $G$-space $(Y,\nu)$. It follows from Furstenberg's Correspondence Principle (see Section \ref{Fur}) that there exists an
ergodic $G$-invariant probability measure $\mu$ on $X$ such that $d^*(A_{x_o}) = \mu(A)$. \\

Towards Theorem \ref{erg1}, we suppose that $\mu(A) > 0$ and that $B \subset Y$ is a Borel set with positive $\nu$-measure 
such that
\begin{equation}
\label{case1}
\nu(A_{x_o}B) < \min(1,d^*(A_{x_o}) + \nu(B)).
\end{equation}
If we define $C = (A_{x_o} B)^c$, then $A_{x_o}^{-1}C \subset B^c$ and 
\[
\nu(A_{x_o}^{-1}C) < d^*(A_{x_o}) + \nu(C).
\]
In Section \ref{sec:actionsets} we show that there exists an ergodic joining $\eta$ of $(X,\mu)$ and $(Y,\nu)$ such that if we write 
\[
A' = A \times Y \qand C' = X \times C
\]
then 
\[
\eta \otimes \eta(G(A' \times C')) < \eta(A') + \eta(C').
\]
In Section \ref{subsect:joincont} and Section \ref{subsec:KM} we further show that there exist
\begin{itemize}
\item a compact group $K$ and a homomorphism $\tau : G \ra K$ with dense image.
\item a closed subgroup $L < K$ and a $G$-factor map $\pi : (X \times Y,\eta) \ra (K/L,m_{K/L})$, 
where $m_{K/L}$ denotes the unique $K$-invariant probability measure on $K/L$.
\item Borel sets $I, J \subset K/L$ with $\nu(A_{x_o}^{-1}C) \geq m_K(I^{-1}J)$, where we have 
identified $I$ and $J$ with their right-$L$-invariant lifts to $K$,
\end{itemize}
such that
\[
A_x \subset \tau^{-1}(I\pi(x,y)^{-1}) \qand C_y \subset \tau^{-1}(J\pi(x,y)^{-1}),
\]
for $\eta$-almost every $(x,y)$. In particular,
\[
m_K(I^{-1}J) < m_K(I) + m_K(J).
\]
Using a classical inverse product set theorem by Kemperman, we conclude that $I^{-1}J$ is 
invariant under an open normal subgroup. It follows that $A_x$ is contained in some proper periodic subset $P$ 
of $G$ for $\mu$-almost every $x$. From this it is not hard to show that there exists a subset $A_o \subset A_{x_o}$ with $d^*(A_o) = d^*(A_{x_o})$
such that $A_o$ is contained in some right translate of $P$. In particular, $A_{x_o}$ is not spread-out, which
finishes the proof of Theorem \ref{erg1} in the case of $A_{x_0}$ being spread-out. 

To prove the second part, we assume that every finite quotient of $G$
is abelian, and we wish to show that not only $A_o$ but the whole of $A_{x_o}$ is contained in some right-translate
of $P$. This is somewhat technical and requires us to use our overshoot inequality \eqref{overshoot} together with some 
classical results of Kneser. \\

The proof of Theorem \ref{erg2} runs along similar lines, but here we end up with Borel sets $I, J \subset K$ such that
\[
m_K(I) = \mu(A) \qand m_K(J) = \nu(B) \qand m_K(I^{-1}J) = m_K(I) + m_K(J) < 1.
\]
A deep fact (see the Appendix in \cite{bj2017}) tells us that since $K$ has a dense \emph{amenable} subgroup, its (possibly trivial) 
identity component $K^o$ must be abelian.
This allows us to use a recent result by the first author which asserts that $K$ admits either $\bT$ or $\bT \rtimes \{-1,1\}$
as a factor in such a way that the sets $I$ and $J$ above coincide, modulo null sets, with pull-backs under the factor map of "intervals".
This shows that for $\mu$-almost every $x \in X$, the set $A_x$ is contained in a Sturmian set with
the same upper Banach density as $A_x$. 

If the $G$-action on $X$ were minimal, or if $\supp(\mu) = X$, then it would follow from
general principles that $A_{x_o}$ is also contained in a Sturmian set with the same upper Banach density. However, these are somewhat degenerate cases. To prove that we can take $x = x_o$ in general requires quite a lot of work (already for abelian $G$), but the necessary arguments are again based on the "overshoot inequality" \eqref{overshoot}.

\addtocontents{toc}{\protect\setcounter{tocdepth}{2}}

\section{Preliminaries}

Throughout this section, let $G$ be a countable group. If $Y$ is a set 
upon which $G$ acts, then we refer to $Y$ as a \emph{G-space}, and if 
$B \subset Y$ and $y \in Y$, then we define the "set of returns" of 
$y$ to $B$ by
\begin{equation}
\label{defreturn}
B_y = \big\{ g \in G \, : \, gy \in B \big\} \subset G.
\end{equation}
We note that for every $B \subset Y$,
\begin{equation}
gB_y = (gB)_y \qand B_{gy} = B_y g^{-1} \quad \textrm{for all $g \in G$ and $y \in Y$}.
\end{equation}
In particular, for every $A \subset G$, we have $AB_y = (AB)_y$ for all $y \in Y$.

\addtocontents{toc}{\protect\setcounter{tocdepth}{1}}
\subsection{Dynamical tools and basic notions}
\label{subsec:basicnot}

\subsubsection{\textbf{Hulls}}
\label{subsec:hull}

Let us denote by $2^G$ the space of all subsets of $G$, endowed with the sequentially compact Tychonoff topology,
with respect to which the set $U = \big\{ A \subset G \, : \, e_G \in A \big\}$ is clopen. We note that the group $G$ acts on $2^G$ 
by homeomorphisms via $g.A = Ag^{-1}$, and using the notation in \eqref{defreturn} above, we have the 
curious-looking identity $U_A = A$. In particular, every subset of $G$ is the set of returns of itself, viewed as an element
in the $G$-space $2^G$, to the set $U$. 

Given $A \subset G$, we shall denote by $X_A$ the closure of the $G$-orbit of 
the point $x_o = A$ in $2^G$. The pair $(X_A,x_o)$ is a pointed $G$-space in the sense 
of Appendix \ref{Fur}, and we refer it as the $G$-\emph{hull} of $A$. It will often be convenient 
to abuse notation and denote by $A$ the clopen set $U \cap X_A$ in $X_A$, so that we can write 
$A_{x_o} = A$.

\subsubsection{\textbf{Borel $G$-spaces and their factors}}
\label{subsec:Gspace}
Let $(Z,\eta)$ be a standard Borel probability measure space with Borel $\sigma$-algebra $\cB_Z$. If it comes equipped
with an action of $G$ by bi-measurable $\eta$-preserving maps, then we say that $(Z,\eta)$ is a 
\emph{Borel $G$-space}. If the only 
$G$-invariant $\cB_Z$-measurable subsets of $Z$ are either $\eta$-null or $\eta$-conull, we say that $\eta$ is \emph{ergodic},
and that $(Z,\eta)$ is an \emph{ergodic Borel $G$-space}. If every finite-index subgroup of $G$ acts ergodically on $(Z,\eta)$
as well, we say that $(Z,\eta)$ is a \emph{totally ergodic Borel $G$-space}.

If $(Z,\eta)$ and $(W,\theta)$ are Borel $G$-spaces, $Z' \subset Z$ and 
$W' \subset W$ are conull $G$-invariant measurable subsets and $\pi : Z' \ra W'$ is a measurable and $G$-equivariant map
such that $\eta(\pi^{-1}(C)) = \theta(C)$ for all $C \in \cB_W$, then we say that $(W,\theta)$ is a \emph{factor Borel $G$-space} of $(Z,\eta)$ and $\pi : (Z,\eta) \ra (W,\theta)$ is a \emph{$G$-factor map}. If $\pi$ in addition is a bijection (which implies that its inverse is measurable as well), we say that $(Z,\eta)$ and $(W,\theta)$
are \emph{isomorphic Borel $G$-spaces}. 

If $\pi : (Z,\eta) \ra (W,\theta)$ is a $G$-factor map, then $\pi^{-1}(\cB_W)$ is a (up to $\eta$-null sets) a $G$-invariant sub-$\sigma$-algebra of $\cB_Z$. Conversely, it is well-known (see for instance Theorem 6.5 in \cite{eiwa2015}), that to every (essentially) $G$-invariant sub-$\sigma$-algebra $\cL \subset \cB_Z$ there correspond 
\begin{enumerate}
\item a factor Borel $G$-space $(W,\theta)$, and
\item a $G$-factor map $\pi : (Z,\eta) \ra (W,\theta)$,
\end{enumerate} 
such that $\cL = \pi^{-1}(\cB_W)$ modulo $\eta$-null sets. We shall refer to any Borel $G$-space $(W,\theta)$ with this property as a \emph{factor Borel $G$-space associated to $\cL$}; they are all isomorphic as Borel $G$-spaces. 

Any Borel $G$-space $(Z,\eta)$ gives rise to a strongly continuous unitary representation of $G$ on the Hilbert space $L^2(Z,\eta)$ via $g \cdot f = f \circ g^{-1}$, which we refer to as the \emph{regular representation} of $(Z,\eta)$; the term \emph{Koopman representation} is often also used in the literature. It is easy to prove that $\eta$ is ergodic if and only if $L^2(Z,\eta)^G \cong \bC$. We say that $(Z,\eta)$ has \emph{discrete spectrum} if $L^2(Z,\eta)$ decomposes into a direct sum of \textsc{finite-dimensional} irreducible representations. 

\subsubsection{\textbf{Discrete spectrum and isometric $G$-spaces}}

Let us now introduce a class of Borel $G$-spaces that will play a key role in this paper. Let 
$(K,\tau)$ be a metrizable compactification of $G$ (see Appendix \ref{sec:gencomp} for the
necessary terminology), and let $L < K$ be a closed subgroup. Then there is an action of $G$ by
continuous maps on the quotient space $K/L$ given by $g.tL = \tau(g)tL$. If $m_{K/L}$ denotes
the unique $K$-invariant probability measure on $K/L$, then $(K/L,m_{K/L})$ is obviously a 
Borel $G$-space (Borel standardness follows from metrizability). We shall refer to $(K,L,\tau)$
as an \emph{isometric $G$-space}; this choice of terminology is standard and comes from 
the fact that there in this setting always exists a $K$-invariant, and thus $G$-invariant, metric on the quotient space $K/L$. 

\begin{lemma}
\label{basic_isomG}
Let $(K,L,\tau)$ be an isometric $G$-space. Then,
\begin{enumerate}
\item $m_{K/L}$ is the unique $G$-invariant Borel probability measure on $K/L$.
\item $G \acts (K/L,m_{K/L})$ is ergodic and has discrete spectrum.
\end{enumerate}
\end{lemma}

\begin{proof}
(i) The stabilizer of any probability measure on $K/L$ is a \emph{closed} subgroup of $K$. Since $\tau(G)$
is dense in $K$, any $G$-invariant probability measure on $K/L$ must be $K$-invariant.

(ii) By continuity of the regular representation of $K$ on $L^2(K/L)$, any $G$-invariant element must be $K$-invariant, and
thus constant. By Peter-Weyl's Theorem, the left-regular representation of $K$ on $L^2(K)$ decomposes into finite-dimensional
irreducible representations, whence the regular representation on $L^2(K/L)$ as well. Each of these finite-dimensional representations is irreducible under the
$G$-action as well (by denseness of $\tau(G)$ in $K$).
\end{proof}

Let us now state a strong converse to Lemma \ref{basic_isomG} (ii) due to Mackey (\cite[Theorem 1]{ma1964}). For \emph{abelian} $G$, this result is more often
referred to as a special instance of the classical Halmos-von Neumann Theorem (see e.g. \cite[Theorem 7.1]{gl2003}).

\begin{proposition}
\label{Mackey}
Every ergodic Borel $G$-space with discrete spectrum is isomorphic as a Borel $G$-space to $(K/L,m_{K/L})$ for some isometric $G$-space $(K,L,\tau)$.
\end{proposition}

Let $(Z,\eta)$ be an ergodic Borel $G$-space, and denote by $\cK$ the smallest $G$-invariant sub-$\sigma$-algebra of $\cB_Z$ with respect to which all elements in $L^2(Z,\eta)$ with finite-dimensional cyclic sub-spaces under the $G$-action are measurable. We note that $\cK = \cB_Z$ if and only if $(Z,\eta)$ has discrete spectrum, as every such finite-dimensional cyclic sub-space decomposes into a direct sum of irreducibles. If $(W,\theta)$ is a factor Borel $G$-space associated to $\cK$, then $G \acts (W,\theta)$ clearly has discrete spectrum, so by Proposition \ref{Mackey} it is isomorphic to $(K/L,m_{K/L})$ for some isometric $G$-space $(K,L,\tau)$. In particular, 
it follows from the discussions in the last subsection that there is a $G$-factor map $\pi : (Z,\eta) \ra (K/L,m_{K/L})$ such that $\pi^{-1}(\cB_{K/L}) = \cK$ modulo $\eta$-null sets. We shall refer to $(K,L,m_{K/L})$ as the \emph{Kronecker-Mackey triple} associated to $(Z,\eta)$, and to both $\cK$ and $(K/L,m_{K/L})$ as the \emph{Kronecker-Mackey factor} of $(Z,\eta)$. The $G$-factor map $\pi$ will be referred to as the \emph{Kronecker-Mackey $G$-factor map}. It is clear from the definition of $\cK$ that any other $G$-invariant sub-$\sigma$-algebra of $\cB_Z$ whose associated factor Borel $G$-space has discrete spectrum is contained in $\cK$. The following lemma will be useful in the next section.

\begin{lemma}
\label{Ginv}
If $\cE_G$ denotes the sub-$\sigma$-algebra of $\cB_Z \otimes \cB_Z$ consisting of $G$-invariant subsets of $Z \times Z$,
then $\cE_G \subset \cK \otimes \cK$.
\end{lemma}

\begin{proof}
We first note that any $G$-invariant function $f \in L^2(Z \times Z,\eta \otimes \eta)$ decomposes as $f_1 + i f_2$, where $f_1$ and $f_2$ are $G$-invariant, and $\overline{f_j(z,z')} = f_j(z',z)$, for $j = 1,2$. It thus suffices to show that any $G$-invariant $f$  with 
$f(z,z') = \overline{f(z',z)}$ is $\cK \otimes \cK$-measurable. To prove this, fix such a (non-zero) $G$-invariant function $f$, and 
consider the (non-zero) operator $T_f : L^2(Z,\eta) \ra L^2(Z,\eta)$ given by
\[
(T_f \phi)(z) = \int_{Z} f(z,z') \phi(z') \, d\eta(z').
\]
It is a well-known classical fact that $T_f$ is self-adjoint and Hilbert-Schmidt. Hence, by the Spectral Theorem for such operators, 
there is an orthonormal basis of eigenfunctions $(\psi_j)$ for $T_f$ with eigenvalues $(\lambda_j)$  
such that $T_f \phi = \sum_j \lambda_j \langle \phi, \psi_j \rangle \psi_j$ for all $\phi \in L^2(Z,\eta)$; in particular, unwrapping this 
identity yields
\begin{equation}
\label{f}
f = \sum_j \lambda_j \psi_j \otimes \overline{\psi_j},
\end{equation}
where convergence holds in the $L^2$-sense. Since $\eta$ is $G$-invariant, $T_f$ is $G$-equivariant, and thus every eigenspace 
of $T_f$ is $G$-invariant. Since $T_f$ is a $G$-invariant operator, each eigenspace is $G$-invariant. By compactness of $T_f$, 
each eigenspace corresponding to a non-zero eigenvalue (such eigenvalues exist since $T_f$ is non-zero) is finite-dimensional.
We conclude that the $G$-cyclic sub-spaces for the corresponding $\psi_j$'s are finite-dimensional, whence $\psi_j$ is $\cK$-measurable. By \eqref{f}, we can conclude that $f$ is $\cK \otimes \cK$-measurable. 
\end{proof}

\subsubsection{\textbf{Shadows}}
\label{subsect:shadows}

Let $(Z,\eta)$ be a Borel $G$-space. Given a sub-$\sigma$-algebra $\cF$ of $\cB_Z$ and a $\cB_Z$-measurable subset $A \subset Z$,
we can consider the conditional expectation $\bE[\chi_A \, | \, \cF]$, pick a $\eta$-almost everywhere defined pointwise realization
of this element in $L^2(Z,\eta)$, and define 
\[
A_{\cF} = \big\{ z \in Z \, : \, \bE[\chi_A \, | \, \cF](z) > 0 \big\}.
\]
We shall refer to any $A_{\cF}$ constructed in this manner as an \emph{$\cF$-shadow of the set $A$}. It is clear that all possible choices of $A_{\cF}$ only differ by $\eta$-null sets, and that $A \subset A_{\cF}$ modulo $\eta$-null sets for all such choices. Moreover, $A_{\cF}$ is $\cF$-measurable by
construction. 

The following lemma will be very useful in the next section. Recall that $\cE_G$ denotes the sub-$\sigma$-algebra of 
$\cB_Z \otimes \cB_Z$ consisting of $G$-invariant subsets of $Z \times Z$.

\begin{lemma}
\label{shadows}
Let $\cF \subset \cB_Z$ be a $G$-invariant sub-$\sigma$-algebra, and suppose that $\cE_G \subset \cF \otimes \cF$. Then, for any $\cB_Z$-measurable sets $A, B \subset Z$, we have, 
modulo $\eta \otimes \eta$-null sets,
\begin{enumerate}
\item $(A \times B)_{\cF \otimes \cF} = A_{\cF} \times B_{\cF}$. 
\item $G(A \times B) = G(A_{\cF} \times B_{\cF})$.
\end{enumerate}
\end{lemma}

\begin{proof}
(i) Obvious; true for any sub-$\sigma$-algebra $\cF$.

(ii) Since $A \subset A_{\cF}$ and $B \subset B_{\cF}$ modulo $\eta$-null sets, it suffices to prove that 
the $G$-invariant set $E = G(A_{\cF} \times B_{\cF}) \setminus G(A \times B)$ is an $\eta \otimes \eta$-null set.
By our assumption on $\cF$, the set $E$ is $\cF \otimes \cF$-measurable, whence
\[
0 = \eta \otimes \eta(E \cap (A \times B)) = \int_{E} \bE[ \chi_{A \times B} \, | \, \cF \otimes \cF ] \, d\eta \otimes \eta
= \int_{E} \bE[ \chi_{A} \, | \cF] \, \bE[ \chi_{B} \, | \, \cF]  \, d\eta \otimes \eta,
\]
and since the integrand on the right is strictly positive on the direct product $A_{\cF} \times B_{\cF}$, we conclude that 
$\eta \otimes \eta(E \cap (A_{\cF} \times B_{\cF})) = 0$. Since $E$ is $G$-invariant and $G$ is countable, 
this implies that $0 = \eta \otimes \eta(E \cap G(A_{\cF} \times B_{\cF})) = \eta \otimes \eta(E)$, which finishes the
proof.
\end{proof}

Combined with Lemma \ref{Ginv}, and the discussion proceeding it, Lemma \ref{shadows} yields the following corollary.

\begin{corollary}
\label{firstreduction}
Let $(Z,\eta)$ be an ergodic Borel $G$-space, and let $(K,L,\tau)$ denote its Kronecker-Mackey triple, and 
$\pi : (Z,\eta) \ra (K/L,m_{K/L})$ the associated $G$-factor map. Then, for all measurable $A, B \subset Z$, 
there are Borel sets $I, J \subset K/L$ such that
\[
A \subset \pi^{-1}(I) \qand B \subset \pi^{-1}(J) 
\]
modulo $\eta$-null sets, and $\eta \otimes \eta(G(A \times B)) = m_{K/L} \otimes m_{K/L}(G(I \times J))$.
\end{corollary}

\begin{proof}
Note that $A \subset A_{\cK}$ and $B \subset B_{\cK}$ modulo $\eta$-null sets, and by Lemma \ref{shadows} (ii) 
we know that, $\eta \otimes \eta(G(A \times B)) = \eta \otimes \eta(G(A_{\cK} \times B_{\cK}))$. Since $A_{\cK}$
and $B_{\cK}$ are $\cK$-measurable, there are Borel sets $I, J \subset K/L$ such that $A_{\cK} = \pi^{-1}(I)$ and
$B_{\cK} = \pi^{-1}(J)$, which finishes the proof.
\end{proof}

\addtocontents{toc}{\protect\setcounter{tocdepth}{2}}

\section{The joining trick}
\label{subsect:joincont}

\begin{adjustwidth*}{1in}{1in}
\emph{We define joinings of Borel $G$-spaces. For compact pointed $G$-spaces, we discuss how one
can use joinings to transfer inclusions of return times at generic points to inclusion of return times 
at base points.}
\end{adjustwidth*}

\vspace{0.2cm}

Let $(X,\mu)$ and $(W,\theta)$ be Borel $G$-spaces. A $G$-invariant Borel probability measure $\xi$ on $X \times W$ such 
that $\mu(A) = \xi(A \times W)$ and $\theta(I) = \xi(X \times I)$ for all Borel sets $A \subset X$ and $I \subset W$ is called a 
\emph{joining} of $(X,\mu)$ and $(W,\theta)$. Note that $\mu \times \theta$ is always a joining. We denote the 
set of all joinings of $(X,\mu)$ and $(W,\theta)$ by $\cJ_G(\mu,\theta)$. The following result is standard (see for 
instance Theorem 6.2 in \cite{gl2003}).

\begin{proposition}
\label{factsjoin0}
If $(X,\mu)$ and $(W,\theta)$ are ergodic, then there are ergodic measures in $\cJ_G(\mu,\theta)$.
\end{proposition}

Let $G$ be a countable amenable group, and let $X$ and $W$ be compact metrizable spaces, equipped with actions 
of $G$ by homeomorphisms, and suppose that there exists a point $x_o \in X$ with a dense $G$-orbit. Let $\mu$ be 
an ergodic $G$-invariant measure on $X$.

\begin{lemma}
\label{conX}
With the notation and assumptions above,
\begin{enumerate}
\item there exists a $\mu$-conull subset $X' \subset X$ such that $\supp(\mu) = \overline{Gx}$ for all $x \in X'$.
\item for every closed $G$-invariant subset $Z \subset X \times W$ whose projection to $X$ contains $\supp(\mu)$,
there exists an ergodic $\xi \in \cP_G(Z)$ which projects to $\mu$. 
\end{enumerate}
In particular, if $\cP_G(W) = \{\theta\}$, then the measure $\xi$ in (ii) is an ergodic joining of $(X,\mu)$ and $(W,\theta)$.
\end{lemma}

\begin{proof}
(i) Since $X$ is compact and metrizable, so is $\supp(\mu)$, and thus there is a countable basis $(U_n)$ for the restricted topology
on $\supp(\mu)$. By ergodicity of $\mu$, we have $\mu(GU_n) = 1$ for all $n$, so $X' = \bigcap_n GU_n$ is $\mu$-conull.
For all $x \in X'$, the $G$-orbit of $x$ meets every $U_n$, and is thus dense in $\supp(\mu)$.

(ii) Fix a closed $G$-invariant set $Z \subset X \times W$, write $p : Z \ra X$ for the projection, and assume that $\supp(\mu) \subset p(Z)$. By a standard use of Hahn-Banach's Theorem, the set of probability measures on $Z$ which project to $\mu$ is non-empty, 
and it is clearly weak*-compact and convex in $\cP(Z)$. Since $G$ is amenable, there is a $G$-fixed point $\xi'$ in this set, and 
since $\mu$ is ergodic, every ergodic component $\xi$ of $\xi'$ will project to $\mu$ as well. 
\end{proof}

Let us also record the following corollary of this lemma; see Appendix \ref{Fur} for terminology concerning pointed $G$-spaces.

\begin{corollary}
\label{super}
Let $G$ be a countable amenable group, $(X,x_o)$ a pointed $G$-space, $\mu$ an ergodic $G$-invariant
probability measure on $X$ and $(K,L,\tau)$ an isometric $G$-space. Then, for every $t \in K/L$ and $\mu$-almost 
every $x \in X$, there exist
\begin{enumerate}
\item an ergodic joining $\xi$ of $(X,\mu)$ and $(K/L,m_{K/L})$ supported on $Z = \overline{G(x,t)}$.
\item a point $t_o \in K/L$ such that $Z \subset Z_o := \overline{G(x_o,t_o)}$.
\end{enumerate}
\end{corollary}

\begin{proof}
(i) is immediate from Lemma \ref{conX} applied to $(W,\theta) = (K/L,m_{K/L})$, using the fact that $m_{K/L}$ is the unique 
$G$-invariant probability measure on $K/L$ (Lemma \ref{basic_isomG}).

(ii) Pick $x \in \supp(\mu)$ and $t \in K/L$, and find a sequence $(g_n)$ in $G$ such that $g_n x_o \ra x$. Choose a 
sub-sequence $(g_{n_k})$ such that $\tau(g_{n_k})$ converges in $K$ to some $k$, and set $t_o = k^{-1} t$. One readily 
checks that $g_{n_k} t_o \ra t$, and thus $(x,t) \in \overline{G(x_o,t_o)}$.
\end{proof}

We now arrive at the punchline of this subsection. In what follows, let $G$ be a countable amenable group, and
\begin{itemize}
\item $(X,x_o)$ pointed $G$-space and $\mu$ an ergodic $G$-invariant probability measure on $X$. 
\item $A \subset X$ an open $\mu$-Jordan measurable subset, i.e. $\mu(\overline{A}) = \mu(A)$.
\item $(K,L,\tau)$ an isometric $G$-space, and $I \subset K/L$ a closed $m_{K/L}$-measurable subset, which
we shall identify with its right-$L$-invariant lift to $K$.
\end{itemize}

Our aim will be to show the following lemma, which roughly asserts that if $A_x \subset I^o_t$ for some $x$ and $t$,
where $I^o$ denotes the interior of $I$, then there is a "big" subset $A_o \subset A_{x_o}$ such that $A_o \subset I_{t_o}$
for some $t_o \in K/L$. This will be the first step in a "bootstrap argument" used in the proofs of Theorem \ref{erg1} and
\ref{erg2}.

\begin{lemma}
\label{containmentlemma}
Suppose that $A_x \subset I^o_t$ for some $x \in X$ such that $\supp(\mu) = \overline{Gx}$ and $t \in K/L$.
Then there exist $t_o \in K/L$ and an extreme invariant mean $\lambda$ on $G$ such that if we set
$A_o = (A \times I^o)_{(x_o,t_o)}$, then
\begin{enumerate}
\item $\lambda(A_o) = \mu(A)$ and $A_o \subset \tau^{-1}(I t_o^{-1})$.
\item for every Borel $G$-space $(Y,\nu)$ and Borel set $C \subset Y$, we have $\nu(A_o^{-1}C) \geq \nu(A_x^{-1}C)$.
\end{enumerate}
\end{lemma}

\begin{proof}
(i) Let $Z = \overline{G(x,t)}$, and use Corollary \ref{super} to produce an ergodic joining $\xi$ of $(X,\mu)$ and $(K/L,m_{K/L})$
supported on $Z$, and a point $t_o \in K/L$ such that $Z \subset \overline{G(x_o,t_o)}$. It follows 
from the inclusion $A_x \subset I_t$ that the \emph{open} set $D = A \times I^c$ satisfies 
$D \cap Z = \emptyset$, whence $\xi(D) = 0$ since $\xi$ is supported on $Z$, and thus $\mu(A) = \xi(A \times I)$. 
Since $\xi$ is an ergodic $G$-invariant measure on $Z_o = \overline{G(x_o,t_o)}$ we can apply Proposition \ref{onto} to the 
pointed $G$-space $(Z_o,z_o)$ where  $z_o = (x_o,t_o)$, and find an extreme invariant mean $\lambda$ such that $\xi = S_{z_o}^*\lambda$. Since $A \times I$ is $\xi$-Jordan measurable, Lemma \ref{somecons} (i) now shows that if 
we write $A_o = (A \times I^o)_{z_o}$, then $\lambda(A_o) = \mu(A)$. The inclusion $A_o \subset \tau^{-1}(I t_o^{-1})$ is immediate.

(ii) Since $A_o \supset U_{z_o}$, where $U = A \times I^o$ is open, the lower bound for every $(Y,\nu)$ and $C \subset Y$ 
follows from the proof of Lemma \ref{symmetry}, by using the fact that $U_z = A_x$ (which is equivalent to the inclusion 
$A_x \subset I_t^o$).
\end{proof}

\subsubsection{\textbf{Removing Jordan measurability}}

In most of our arguments, assuming that the set $A \subset X$ is $\mu$-Jordan measurable is rather harmless; in fact, in most
of our applications, $A$ will be clopen. However, at one subtle point in the proof of our main density results, it will be useful 
to instead refer to the following weak cousin to Lemma \ref{containmentlemma}. 

\begin{lemma}
\label{convtolemma}
Suppose that $U \subset X$ is open, and there is a Borel set $Q \subset U$ with positive $\mu$-measure which is invariant
under a finite-index normal subgroup $G_o < G$. Then there is a thick set $T \subset G$ and a non-empty $G_o$-invariant set $Q_o \subset G$
such that $Q_o \cap T \subset U_{x_o}$.
\end{lemma}

\begin{proof}
One readily checks that the map $\sigma : X \ra 2^{G/G_o}$ given by $x \mapsto Q_x$ is Borel and $G$-equivariant, and $Q = \sigma^{-1}(V)$,
where $V = \{ D \subset G/G_o \, : \, G_o \in D \big\}$ is clopen. In particular, we have $V_{\sigma(x)} \subset U_x$ for all $x \in X$,
and thus $\xi(U^c \times V) = 0$ for the graph joining $\xi = (\textrm{id} \times \sigma)_*\mu$. Just as in Corollary \ref{super},
we can utilize the ergodicity of $\xi$ to find $t_o \in K/L$ such that $\overline{G(x_o, t_o)} \supset \supp(\xi)$. We then use Proposition \ref{onto} to find an invariant mean $\lambda$ such that $S_{(x_o,t_o)}^*\lambda = \xi$, where $S_{(x_o,t_o)}^*$ is defined in Subsection \ref{subsec:pointedG}. Since $U^c \times V$ is closed, Lemma \ref{somecons} tells us that 
$\lambda(U^c_{x_o} \cap V_{t_o}) \leq \eta(U^c \times V) = 0$, and thus $T = (U^c_{x_o} \cap V_{t_o})^c$ is thick by Lemma \ref{thicksyn}. The set $Q_o = V_{t_o}$ is clearly $G_o$-invariant, and one readily checks that $Q_o \cap T \subset U_{x_o}$.
\end{proof}

\addtocontents{toc}{\protect\setcounter{tocdepth}{1}}
\subsection{Extra features of the joining trick if the action is minimal (optional)}
\addtocontents{toc}{\protect\setcounter{tocdepth}{2}}

If one were to impose on the set $A$ in Theorem \ref{erg1} or Theorem \ref{erg2} the additional (somewhat unnatural) assumption that its $G$-hull
is a minimal $G$-space (or at least that there exists a $G$-invariant measure on the hull of full support), then many arguments in the coming sections would become significantly shorter and less technical, in view of the following
stronger version of Lemma \ref{containmentlemma}. We retain the notation introduced in the previous subsection.

\begin{lemma}
\label{containmentlemmaminimal}
Suppose that $G \acts X$ is minimal, $A \subset X$ is open and $I \subset K/L$ is closed. 
If $A_x \subset I_t$ for some $(x,t) \in X \times K/L$, then there exists $t_o \in K/L$ such that
$A_{x_o} \subset I_{t_o}$.
\end{lemma}

\begin{proof}
Let $Z = \overline{G(x,t)} \subset X \times K/L$. Since $X$
is minimal, we conclude that $Z$ projects onto $X$, and thus there is at least one $t_o \in K/L$ such that $(x_o,t_o) \in Z$. 
We further note that since $A \times I^c$ is open, and $(A \times I^c)_{(x,t)} = \emptyset$, we have $(A \times I^c) \cap Z = \emptyset$,
and thus in particular $A_{x_o} \cap I^c_{t_o} = \emptyset$, whence $A_{x_o} \subset I_{t_o}$.
\end{proof}

\section{How to take differences in Borel $G$-spaces}
\label{subsec:KM}

\begin{adjustwidth*}{1in}{1in}
\emph{We discuss a way to 
associate to pairs of Borel sets in an ergodic Borel $G$-space their
"difference set" in the associated Kronecker-Mackey factor.}
\end{adjustwidth*}

\vspace{0.2cm}

Let $G$ be a countable group, and $(Z,\eta)$ an ergodic Borel $G$-space. Let $A, B \subset Z$ be 
Borel sets. If $(K,L,\tau)$ denotes the Kronecker-Mackey triple associated to $(Z,\eta)$ (see previous
section for definitions), and $\pi : (Z,\eta) \ra (K/L,m_{K/L})$ the corresponding $G$-factor map, then 
Corollary \ref{firstreduction} provides us with Borel sets $I, J \subset K/L$ such that
$A \subset \pi^{-1}(I)$ and $B \subset \pi^{-1}(J)$ and
\[
\eta \otimes \eta(G(A \times B)) = m_{K/L} \otimes m_{K/L}(G(I \times J)).
\]
We may identify $I$ and $J$ with their right-$L$-invariant lifts to $K$ under the canonical map $K \mapsto K/L$,
after which we can write this identity as $\eta \otimes \eta(G(A \times B)) = m_{K} \otimes m_{K}(G(I \times J))$.
It is now tempting to argue as follows: Since $\tau(G)$ is dense in $K$, we should be able to replace $G(I \times J)$ with 
$K(I \times J)$ without increasing the $\eta \otimes \eta$-measure; the latter set is the pull-back of the set $I^{-1}J$
under the multiplication map $(x,y) \mapsto x^{-1}y$, and since $m_K \otimes m_K$ is mapped to $m_K$ under this
multiplication map, we have $\eta \otimes \eta(G(A \times B)) = m_{K}(I^{-1}J)$. Unfortunately, already the first
line of the argument above fails; replacing $I$ and $J$ with $I \cup N_I$ and $J \cup N_J$ where $N_I$ and $N_J$ 
are $m_K$-null sets such that the difference set $N_I^{-1}N_J$ has positive $m_K$-measure shows that additional 
arguments are required. Fortunately for us, upon passing to conull subsets of $I$ and $J$ in the first step (which will 
not affect the measure of $G(I \times J)$), the rest of the argument runs as before. The exact correction can be stated
as follows. 

\begin{proposition}
\label{IJ}
If $I, J \subset K$ are Borel sets, then there exist conull subsets $I' \subset I$ and $J' \subset J$ such that
\[
m_K \otimes m_K(G(I' \times J')) = m_K \otimes m_K(K(I' \times J')).
\]
Furthermore, if $I$ and $J$ are right-invariant under a subgroup $L < K$, then so are $I'$ and $J'$.
\end{proposition}

\subsection{Difference arithmetics for shadows}

Combining Proposition \ref{IJ} with the rest of the argument above, as well as with Corollary \ref{firstreduction}, yields
the following corollary which will play a key role in this paper. We stress that this result is new already for actions of
$G = (\bZ,+)$.

\begin{corollary}
\label{KM-products of sets0}
Let $(Z,\eta)$ be an ergodic Borel $G$-space, let $(K,L,\tau)$ denote its Kronecker-Mackey triple, and 
let $\pi : (Z,\eta) \ra (K/L,m_{K/L})$ denote the associated $G$-factor map. Then, for all measurable $A, B \subset Z$, 
there are Borel sets $I, J \subset K/L$ such that
\[
A \subset \pi^{-1}(I) \qand B \subset \pi^{-1}(J) 
\]
modulo $\eta$-null sets, and $\eta \otimes \eta(G(A \times B)) = m_{K}(I^{-1}J)$, where we have identified $I$
and $J$ with their right-$L$-invariant lifts to $K$ under the canonical map $K \mapsto K/L$.
\end{corollary}

The rest of this section will be devoted to the proof of Proposition \ref{IJ}. It will be useful to adopt a slightly more
general perspective. In what follows, let $H$ be a compact metrizable group, $M < H$ a closed subgroup
and $\Gamma < M$ a dense countable subgroup. Soon enough, we shall apply our results below to the setting:
\begin{equation}
\label{grptriple}
H = K \times K \qand M = \big\{ (k,k) \, : \, k \in K \big\} \qand \Gamma = \big\{(\tau(g),\tau(g)) \, : \, g \in G \big\}.
\end{equation}
We say that a decreasing sequence $(B_j)$ of closed sets in $H$ with positive
$m_H$-measures is 
\emph{Dirac} if their intersection equals $\{e_H\}$, and we say that a Borel set $D \subset H$ is \emph{balanced
with respect to $(B_j)$} if 
\[
\lim_j \frac{m_H(D \cap sB_j)}{m_H(B_j)} = 1, \quad \textrm{for all $s \in D$}.
\]
It is easy to see that every compact and metrizable group admits a Dirac sequence $(B_j)$, and given any such 
sequence, we can form $\rho_j = \frac{\chi_{B_j}}{m_H(B_j)}$ in $L^1(H)$. It is quite standard to prove that for 
every Borel set $D \subset H$, we have $\rho_j * \chi_D \ra \chi_D$ in the $L^2$-norm, and thus, upon passing
to a sub-sequence, $m_H$-almost everywhere. Clearly, if $D$ is \emph{right}-invariant under some subgroup $Q$ of $H$, 
then the set on which this sub-sequence converges is again right-invariant under $Q$. Unwrapping this, we conclude:
\begin{lemma}
\label{bal}
If $(B_j)$ is Dirac and $D \subset H$ is Borel with positive $m_H$-measure, then there exists a conull subset 
$D' \subset D$ and a sub-sequence $(B_{j_k})$ such that $D'$ is balanced with respect to $(B_{j_k})$. If $Q$
is a subgroup of $H$, and $D$ is right-$Q$-invariant, then so is $D'$.
\end{lemma}

Let us now assume that $D \subset H$ is a Borel set with positive $m_H$-measure which is balanced with respect to 
some Dirac sequence $(B_j)$. We claim that $m_H(\Gamma D) = m_H(M D)$. To prove this, we argue by contradiction, 
and define $C = MD \setminus \Gamma D$, and assume that $C$ has positive $m_H$-measure. Then for every $j$, 
the function $f_j(s) = m_H(C \cap sB_j)$ is continuous and left $\Gamma$-invariant, whence left $M$-invariant as well.
Furthermore, fix $0 < \eps < 1/2$, and use the lemma above to produce $C' \subset C$ with $m_H(C') = m_H(C)$ which is balanced
with respect to some sub-sequence $(B_{j_k})$. Fix $s \in C'$ and write $s = md$ for some $m \in M$ and $d \in D$. 
Then, since both $C$ and $D$ are balanced with respect to $(B_{j_k})$, we have for large $k$,
\[
f_{j_k}(d) = m_H(C \cap dB_{j_k}) = f_{j_k}(md) = m_H(C \cap sB_{j_k}) \geq (1-\eps) m_H(B_{j_k}),
\]
and
\[
m_H(D \cap dB_{j_k}) \geq (1-\eps) m_{H}(B_{j_k}),
\]
and thus $m_H(C \cap D \cap dB_{j_k}) \geq (1-2\eps) m_H(B_{j_k}) > 0$. In particular, $C \cap D \neq \emptyset$, which
is a contradiction, and thus $m_H(C) = 0$.

Let us now apply all of this to prove Proposition \ref{IJ}, using $H, M$ and $\Gamma$ from \eqref{grptriple}. Let $I, J \subset K$
be Borel sets, and fix a Dirac sequence $(B_j)$. By using Lemma \ref{bal} twice, we can produce a sub-sequence $(B_{j_k})$
and $I' \subset I$ and $J' \subset J$ with $m_K(I') = m_K(I)$ and $m_K(J') = m_K(J)$ which are both balanced with respect to 
$(B_{j_k})$. Note that they can both be chosen to be right-$L$-invariant. We conclude that $D = I' \times J'$ is balanced with 
respect to $B_{j_k} \times B_{j_k}$ in $H = K \times K$, so the argument above tells us that 
\[
m_K \otimes m_K(K(I' \times J')) = m_H(M D) = m_H(\Gamma D) = m_K \otimes m_K(G(I' \times J')).
\]

\section{Action sets versus product sets in compact groups}
\label{sec:actionsets}

\begin{adjustwidth*}{1in}{1in}
\emph{We show that action sets for an ergodic action of a countable group $G$ 
can be "nicely shadowed" by product sets in a group compactification of $G$. 
In certain situations, when a priori upper bounds are imposed on the action sets, topological 
regularity for the involved sets in the group compactification can be deduced.}
\end{adjustwidth*}

\vspace{0.2cm}

This is a long and somewhat technical section, which we partition into three main sub-sections. The same notation is 
kept throughout the section, but with each new sub-section, additional assumptions on the basic objects will be imposed. 

\subsection{A correspondence principle for action sets}

Let $G$ be a countable, not necessarily amenable, group. Throughout this section, our key players will be:
\begin{itemize}
\item a pointed $G$-space $(X,x_o)$ and an ergodic $\mu \in \cP_G(X)$
\item a non-empty open set $A \subset X$.
\item an ergodic Borel $G$-space $(Y,\nu)$ and a Borel set $C \subset Y$ with positive $\nu$-measure.
\end{itemize}
With this notation understood, we define $A', C' \subset X \times Y$ by
\begin{equation}
\label{defAprimCprim}
A' = A \times Y \qand C' = X \times C.
\end{equation}
Our first goal will be to establish the following theorem, which is one of the principal
building blocks in the proofs of our main results.  The theorem admits many immediate, yet 
interesting corollaries. These are stated in Section \ref{sec:mainresults}. 

\begin{theorem}
\label{corr1}
For every ergodic joining $\eta \in \cJ_G(\mu,\nu)$, there exist
\begin{enumerate}
\item a metrizable compactification $(K,\tau)$ of $G$,
\item a closed subgroup $L < K$ and a $G$-factor map $\pi : (X \times Y,\eta) \ra (K/L,m_{K/L})$,
\item Borel sets $I, J \subset K/L$,
\end{enumerate}
such that 
\begin{equation}
\label{inclusions}
A' \subset \pi^{-1}(I) \qand C' \subset \pi^{-1}(J), \quad \textrm{modulo $\eta$-null sets},
\end{equation}
and 
\begin{equation}
\label{lowerbnd}
\nu(A_{x_o}^{-1}C) \geq \eta \otimes \eta(G(A' \times C')) = m_{K}(I^{-1}J), 
\end{equation}
where we have identified $I$ and $J$ with their right $L$-invariant lifts to $K$ under the canonical quotient map $K \ra K/L$. In particular, we have $m_K(I) \geq \mu(A)$ and $m_K(J) \geq \nu(C)$.
\end{theorem}

\begin{remark}
The two Borel $G$-spaces $(X,\mu)$ and $(Y,\nu)$ could be quite different. The choice of an ergodic joining between these two spaces allows us to put them on an equal footing. The price we pay is that we have consider $G$-factors of the bigger space $(X \times Y,\eta)$ instead of $G$-factors of $(X,\mu)$ and $(Y,\nu)$ respectively.
\end{remark}

Let us start the proof of Theorem \ref{corr1} by picking an ergodic $\eta \in \cJ_G(\mu,\nu)$ once and for all. 
The following lemma symmetrizes the roles of $A$ and $C$ so that the results of 
Section \ref{subsec:KM} can be applied. It is the only place in the proof of Theorem \ref{corr1} where the assumption 
that $A$ is open is used. 

\begin{lemma}
\label{symmetry}
With the notation above, 
\begin{enumerate}
\item $\nu(A_{x_o}^{-1}C) \geq \eta \otimes \eta(G(A' \times C'))$.
\item there is a $\mu$-conull subset $X' \subset X$ such that $\nu(A_x^{-1}C) = \eta \otimes \eta(G(A' \times C'))$
for all $x \in X'$.
\end{enumerate}
\end{lemma}

\begin{proof}
First note that the continuous map $\sigma : (X \times Y)^2 \ra X \times Y$ given by 
$((x_1,y_1), (x_2,y_2)) \mapsto (x_1,y_2)$ is $G$-equivariant, and satisfies
\[
\sigma_*(\eta \otimes \eta) = \mu \otimes \nu \qand \sigma^{-1}(A \times C) = A' \times C'.
\]
Hence, $\eta \otimes \eta(G(A' \times C')) = \mu \otimes \nu(G(A \times C))$. Secondly,
\[
\mu \otimes \nu(G(A \times C)) = \int_X \Big( \int_Y  \chi_{G(A \times C)}(x,y) \, d\nu(y) \Big) \, d\mu(x)  
= \int_X \nu(A_x^{-1}C) \, d\mu(x).
\]
Since $\nu$ is $G$-invariant and $A_{gx} = A_x g^{-1}$ for all $g \in G$ and $x \in X$, the measurable function $x \mapsto \nu(A_x^{-1}C)$ is $G$-invariant. Hence, by ergodicity of $\mu$, there exists
a $\mu$-conull subset $X' \subset X$ on which this function equals its $\mu$-integral, which in this case is 
$\eta \otimes \eta(G(A' \times C'))$, showing (ii). 

To prove (i), it suffices to 
establish the lower bound $\nu(A_{x_o}^{-1}C) \geq \nu(A_x^{-1}C)$ for all $x \in X$, as (i) then follows upon integration
against $\mu$. To show this lower bound, let us 
fix $\eps > 0$ and $x \in X$, and assume, without loss of generality, that $A_x$ is non-empty. By $\sigma$-additivity 
of the measure $\nu$, we can find a finite subset $F \subset A_x$ such that $\nu(F^{-1}C) \geq \nu(A_x^{-1}C) - \eps$. 

Since the set $A$ is assumed to be open, the non-empty set $A_F = \{ z \in X \, : \, F \subset A_z \big\}$
is open as well. Since $x_o$ has a dense $G$-orbit in $X$, there exists at least one $g \in G$ such that $g.x_o \in A_F$, whence $Fg \subset A_{x_o}$. Finally, since $\nu$ is $G$-invariant, we conclude that 
$\nu(A_{x_o}^{-1}C) \geq \nu(F^{-1}C) \geq \nu(A_x^{-1}C) - \eps$. Since $\eps > 0$ and $x \in X$ were arbitrary, 
we are done.
\end{proof}

\subsubsection{Proof of Theorem \ref{corr1}}

Let $(K,L,\tau)$ be the Kronecker-Mackey factor of 
$(X \times Y,\eta)$ and denote by $\pi$ the corresponding $G$-factor map $(X \times Y,\eta) \ra (K/L,m_{K/L})$.
By Corollary \ref{KM-products of sets0} applied to $(X \times Y,\eta)$, we can find $I, J \subset K/L$ such that
\[
A' \subset \pi^{-1}(I) \qand C' \subset \pi^{-1}(J), \quad \textrm{modulo $\eta$-null sets},
\]
and, $\eta \otimes \eta(G(A' \times C')) = m_K(I^{-1}J)$, so in combination with the lemma above, 
\[
\nu(A_{x_o}^{-1}C) \geq \eta \otimes \eta(G(A' \times C')) = m_K(I^{-1}J),
\]
where $I$ and $J$ have been identified with their right $L$-invariant lifts to $K$ under the canonical quotient map 
$K \ra K/L$.

\subsection{Forcing regularity from small doubling}

The guiding question in this subsection is:
\vspace{0.2cm}
\begin{center}
\emph{What can be said about the Borel sets $I, J \subset K$ that we end up with in Theorem \ref{corr1}?}
\end{center}  
\vspace{0.2cm}
In general, the answer is "very little". However, as we will see below, if one \emph{assumes} certain \emph{a priori} upper bounds on 
$\nu(A_{x_o}^{-1}C)$ in terms of $d^*(A_{x_o})$ and $\nu(C)$, then powerful tools from the research field of product set theory in groups become available, and will force 
$I$ and $J$ to be "nicely contained" in highly regular sets. This regularity will allow us in the next subsection to utilize the arguments 
from Section \ref{subsect:joincont} ("The joining trick")  in order to establish relations between the sets $A_{x_o}$ and 
$I_{t_o}$ for a certain point $t_o \in K/L$. However, before we state our second main result of this section, we will need some technical definitions.

\begin{definition}
\label{def:reduction}
Let $K$ be a compact metrizable group with Haar probability measure $m_K$, and let $I, J \subset K$
be Borel sets. We say that
\begin{enumerate}
\item $(I,J)$ \emph{reduces} to a pair $(I_o,J_o)$ of Borel sets in a factor group $M$ of $K$ if 
\begin{equation}
\label{reduction}
I \subset p^{-1}(I_o) \qand J \subset p^{-1}(J_o) \qand m_K(I^{-1}J) = m_M(I_o^{-1}J_o),
\end{equation}
where $p : K \ra M$ denotes the factor map. 
\item $(I,J)$ is \emph{left-balanced} if the inclusion $s^{-1}J \subset I^{-1}J$ implies that $s \in I$.
\item $I \subset K$ is \emph{periodic} if it is invariant under an open normal subgroup.
\item $I \subset K$ is \emph{sub-periodic} if 
there exists a conull Borel subset $I' \subset I$ and a normal open subgroup 
$U$ of $K$ such that $I'U \neq K$. 
\item $I \subset K$ is \emph{Sturmian} if either 
\begin{enumerate}
\item $K = \bT$ and $I$ is a closed interval.
\item $K = \bT \rtimes \{-1,1\}$ and $I = (I' \rtimes \{-1,1\})k$ for some symmetric closed interval $I' \subset \bT$, 
and $k \in K$.
\end{enumerate}
\item $(I,J)$ is a \emph{Sturmian pair} if both $I$ and $J$ are Sturmian in the sense of (a) or (b) simultaneously.
\end{enumerate}
\end{definition}

Throughout the rest of the section we identify $I$ and $J$ from Theorem \ref{corr1} with their right $L$-invariant lifts 
to $K$ under the canonical quotient map from $K \ra K/L$. The important role of the theorem below will hopefully 
become clear in the next sub-section.

\begin{theorem}
\label{corr2}
With the notation and assumptions of Theorem \ref{corr1}, we have
\begin{enumerate}
\item $\nu(A_{x_o}^{-1}C) < 1$ $\implies$ $\overline{I} \neq K$. 
\item $\nu(A_{x_o}^{-1}C) < \min(1,\mu(A) + \nu(C))$ $\implies$ $(I,J)$ reduces to $(I_o,J_o)$ in a 
\emph{finite} quotient group $M$ of $K$. In particular, $I$ is contained in a proper periodic subset of $K$. 
Furthermore, if $M$ is abelian, then we can choose $(I_o,J_o)$ to be left-balanced
and satisfy
\begin{equation}
\label{knesercorr}
m_M(I_o^{-1}J_o) = m_M(I_o) + m_M(J_o) - m_M(\{e_M\}).
\end{equation}
\item $G$ is amenable and $\nu(A_{x_o}^{-1}C) = \mu(A) + \nu(C) < 1$ $\implies$ Either $I$ or $J$ is a sub-periodic set, or the pair
$(I,J)$ reduces to a Sturmian pair. In the latter case, we also have $m_K(I) = \mu(A)$ and $m_K(J) = \nu(C)$.
\end{enumerate}
\end{theorem}

We see that in each of the sub-cases in the theorem above, an a priori bound on $\nu(A_{x_o}^{-1}C)$ forces some regularity on 
$I$ and $J$. Theorem \ref{corr1} tells us that such bounds immediately imply analogous bounds on $m_K(I^{-1}J)$. By using 
the following series of results of Kemperman \cite{ke1964}, Kneser \cite{kn1956} and the first author \cite{bj2017}, the proof of 
Theorem \ref{corr2} will be swift. 

\begin{theorem}
\label{potpurri}
Let $K$ be a compact metrizable group with identity component $K^o$, and let $I, J \subset K$ be Borel sets with positive measures.
\begin{enumerate}
\item $m_K(I^{-1}J) < 1$ $\implies$ $\overline{I} \neq K$.
\item $m_K(I^{-1}J) < \min(1,m_K(I) + m_K(J))$ $\implies$ $(I,J)$ reduces to $(I_o,J_o)$ in a 
\emph{finite} quotient group $M$ of $K$. In particular, $I$ is contained in a proper periodic subset of $K$. 
Furthermore, If $M$ is abelian, then we can choose $(I_o,J_o)$ to be left-balanced
and satisfy
\[
m_M(I_o^{-1}J_o) = m_M(I_o) + m_M(J_o) - m_M(\{e_M\}).
\]
If we assume that every finite quotient of $G$ is \emph{abelian}, then $M$ is abelian.
\item $K^o$ is abelian, and $m_K(I^{-1}J)  = m_K(I) + m_K(J) < 1$ $\implies$ Either $I$ or $J$ is a 
sub-periodic set, or the pair $(I,J)$ reduces to a Sturmian pair. 
\end{enumerate}
\end{theorem}

\begin{remark}
We encourage the reader to verify that if $(I,J)$ is Sturmian then neither $I$ nor $J$ is sub-periodic and 
$m_K(I^{-1}J) = \min(1,m_K(I) + m_K(J))$. 
\end{remark}

Concerning the exact credits in Theorem \ref{potpurri}: The first assertion in (ii) is due to Kemperman (Theorem 1, \cite{ke1964}), 
while the second assertion in (ii) is due to Kneser (Satz 1, \cite{kn1956}), modulo the comment about left-balance; this follows from
the simple observation: 
\begin{lemma}
\label{left-balance}
Let $M$ be a compact group, and let $I_o, J_o \subset M$ be closed sets. Then there
exists a closed set $I_1 \subset M$ such that 
\[
I_o \subset I_1 \qand I_o^{-1}J_o = I_1^{-1} J_o \qand \textrm{$(I_1,J_o)$ is left-balanced}.
\] 
\end{lemma}
\begin{proof}
Define $I_1^{-1} := \bigcap_{y \in J_o} I_o^{-1} J_o y^{-1}$, and verify the conditions above. 
\end{proof}
Going back to the credits in Theorem \ref{potpurri}: In the case when $K$ is connected, and thus $K^o = K$, 
(iii) is due to Kneser (Satz 4, \cite{kn1956}). Note that in this case, sub-periodic subsets do not exist. The general 
case of (iii) is due to the first author (Theorem 1.8, \cite{bj2017}). 

Finally, (i) should be attributed to Weil \cite{We} (based on an earlier observation by Steinhaus), 
although this exact form is not stated there. However, it is not hard to deduce "our" version: Note that if $I$ is dense 
in $K$, but $D = (I^{-1}J)^c$ has positive Haar measure, then the dense set $I^{-1}$ would intersect the product set 
$JD^{-1}$ trivially; however, Weil shows that $JD^{-1}$ always has non-empty interior, which leads to a contradiction. 
Hence $m_K(I^{-1}J) = m_K(D^c) = 1$.

\subsubsection{Proof of Theorem \ref{corr2}}

Before we begin, recall from Theorem \ref{corr1} that 
\[
m_K(I) \geq \mu(A) \qand m_K(J) \geq \nu(C) \qand \nu(A_{x_o}^{-1}C) \geq m_K(I^{-1}J),
\]
and $(K,\tau)$ is a metrizable compactification of $G$. 

Hence, if $\nu(A_{x_o}^{-1}C) < 1$, then $m_K(I^{-1}J) < 1$ as well, which by Theorem \ref{potpurri} settles (i). If 
$\nu(A_{x_o}^{-1}C) < \min(1,\mu(A) + \nu(C))$, then 
$m_K(I^{-1}J) < \min(1,m_K(I) + m_K(J))$ as well, which by Theorem \ref{potpurri} settles (ii). Concerning (iii), we note
that if $\nu(A_{x_o}^{-1}C) = \mu(A) + \nu(C) < 1$, then 
\[
m_K(I^{-1}J) \leq \nu(A_{x_o}^{-1}C) = \mu(A) + \nu(C) \leq \min(1,m_K(I) + m_K(J)).
\]
If any of these inequalities is strict, then Theorem \ref{potpurri} (ii) implies that $I$ is contained in a proper periodic
subset of $K$, which in particular means that $I$ is sub-periodic. Hence, if we assume that neither $I$ nor $J$ is 
sub-periodic, then we must have $m_K(I^{-1}J) = m_K(I) + m_K(J) < 1$. Since $G$ is amenable and $(K,\tau)$ is 
a group compactification of $G$, Lemma \ref{compactify} (i) 
guarantees that $K^o$ is abelian, so Theorem \ref{potpurri} (iii) now tells us that $(I,J)$ reduces to a Sturmian pair. 

\subsection{Proving containment using the joining trick}

Let us now take a closer look at what the inclusions \eqref{inclusions} and Theorem \ref{corr2} together imply for the 
set $A_{x_o} \subset G$. To get interesting results, it is necessary to assume from now on that 
\vspace{0.2cm}
\begin{center}
\emph{$A$ is not only open, but also $\mu$-Jordan measurable.}
\end{center}  
\vspace{0.2cm}
From the inclusions \eqref{inclusions}, we know that $A \times Y \subset \pi^{-1}(I)$ modulo $\eta$-null sets. We conclude 
that for $\eta$-almost every $(x,y) \in X \times Y$, we have 
\begin{equation}
\label{coninc1}
A_x = \big(A \times Y)_{(x,y)} \subset \pi^{-1}(I)_{(x,y)} = \tau^{-1}(I \pi(x,y)^{-1}),
\end{equation}
where we have identified $I$ with its right-$L$-invariant lift to $K$. Of course, we also have
\begin{equation}
\label{coninc2}
C_y = (X \times C)_{(x,y)} \subset \tau^{-1}(J\pi(x,y)^{-1}), \quad \textrm{for $\eta$-almost every $(x,y)$}.
\end{equation}
 
As a warm-up for the things that will come, let us first assume that $I$ is sub-periodic in $K$, that is to say, there is a conull subset $I' \subset I$ and an open subgroup $U < K$ such that $U I' \neq K$. By passing to finite-index subgroups, we may without loss of generality assume that $U$ is normal in $K$.  Let $p_o$ denote the canonical quotient map from $K$ to the finite group $K/U$, and 
let $I_o$ denote the image of $I'$ under $p_o$. Then $I_o \neq K/U$, and $A_x \subset \tau_o^{-1}(I_o \pi(x,y)^{-1})$ for $\eta$-almost every $(x,y)$, where $\tau_o = p_o \circ \tau$. Since $K/U$ is finite, $I_o$ is definitely closed and $m_{K/U}$-Jordan measurable. 

Lemma \ref{containmentlemma} now shows that there exists a subset $A_o \subset A_{x_o}$ and an extreme invariant mean $\lambda$ on $G$ such that $\lambda(A_o) = \lambda(A_{x_o})$ and $A_o \subset P = \tau_o^{-1}(I_o t_o^{-1})$, where $P$ is a proper subset, invariant under the finite-index normal subgroup $G_o = \tau^{-1}(U)$. Furthermore, if we in addition assume that the $G$-action on 
$X$ is minimal (or if $\mu$ has full support), then Lemma \ref{containmentlemmaminimal} tells us that we in fact have the stronger inclusion $A_{x_o} \subset P$.

If $J$ is sub-periodic, then we can argue along the same lines as above; if $J' \subset J$ is conull and $U < K$ is an open 
subgroup such that $UJ' \neq K$, then $G_o C_y \subset \tau^{-1}(UJ' \pi(x,y)^{-1}) \neq G$, with $G_o = \tau^{-1}(U)$,
for $\nu$-generic $y \in Y$, and thus $G_o C$ cannot be $\nu$-conull, as we would then have $G_o C_y = G$ for $\nu$-almost
every $y$. \\

We summarize these observations in the following proposition.

\begin{proposition}
\label{mainpropcont1}
With the notation above:
\begin{enumerate}
\item If $I$ is sub-periodic, then there exist a finite-index subgroup $G_o < G$, an extreme invariant mean $\lambda$ on $G$,
and a subset $A_o \subset A_{x_o}$ with $\lambda(A_o) = \mu(A)$ such that $G_o A_o \neq G$. In particular, if $A_{x_o}$ is spread-out, 
then $I$ cannot be sub-periodic. If we in addition assume that 
$G \acts X$ is a minimal action, then $G_o A_{x_o} \neq G$. 
\item If $J$ is sub-periodic, then there is a finite-index subgroup $G_o < G$ such that $\nu(G_o C) < 1$. In particular, if 
$G \acts (Y,\nu)$ is totally ergodic, then $J$ cannot be sub-periodic.
\end{enumerate}
\end{proposition}

\subsubsection{\textbf{The overshoot bound}}

Let us now assume that the pair $(I,J)$ of Borel sets in $K$ from Theorem \ref{corr1} reduces to a pair $(I_o,J_o)$ in a 
quotient group $M$ of $K$ under the quotient map $p : K \ra M$, and let us further assume that $I_o$ and $J_o$ are both 
closed and $m_M$-Jordan measurable. Then, since $I \subset p^{-1}(I_o)$, the inclusions in \eqref{inclusions} imply that
\[
A_x = (A \times Y)_{(x,y)} \subset \tau^{-1}(p^{-1}(I_o)\pi(x,y)^{-1}) = \tau_o^{-1}(I_o \pi_p(x,y)^{-1}), 
\]
for $\eta$-almost every $(x,y)$, where $\tau_o = p \circ \tau$ and $\pi_p = p \circ \pi$.

Upon passing to a further $\eta$-conull subset
we can even ensure that $A_x \subset \tau_o^{-1}(I^o_o t^{-1})$ for $\eta$-almost every $(x,y)$, where $t = \pi_p(x,y)$. 
In particular, the conditions of Lemma \ref{containmentlemma} are satisfied, so we conclude that there exist $t_o \in M$,
an extreme invariant mean $\lambda$ on $G$ and $A_o \subset A_{x_o}$ such that
\[
\lambda(A_o) = \mu(A) \qand A_o \subset \tau_o^{-1}(I_o t_o^{-1}).
\]
and $\nu(A_o^{-1}C) \geq \nu(A_x^{-1}C)$. We can further choose $x$ so that it belongs to the set $X'$ in Lemma \ref{symmetry},
ensuring that 
\[
\nu(A_{x}^{-1}C) = \eta \otimes \eta(G(A \times C)) = m_K(I^{-1}J) = m_M(I_o^{-1}J_o),
\] 
where $A' = A \times Y$ and $C' = X \times C$. In particular, if $M$ is finite, then $A_o$ is contained in a proper periodic subset. If $G \acts X$ is minimal, Lemma \ref{containmentlemmaminimal} tells us that we can choose $t_o \in M$ such that 
$A_{x_o} \subset \tau_o^{-1}(I_o t_o^{-1})$. Without the assumption of minimality, this inclusion might not hold, and 
we have to take a different route. Crucial to this alternative is an overshoot-inequality which we now formulate. 

\begin{lemma}
\label{lemma_overshoot}
For all $s \in A_{x_o} \setminus A_o$, 
\begin{equation}
\label{overshoot}
\nu(A_{x_o}^{-1}C) - \nu(C) \geq m_M(I_o^{-1}J_o) - m_M(J_o) + m_M( \tau_o(s)^{-1} J_o \setminus t_o I_o^{-1}J_o).
\end{equation}
\end{lemma}

\begin{proof}
Pick $s \in A_{x_o} \setminus A_o$, and note that, since $\eta$ is a joining of $(X,\mu)$ and $(Y,\nu)$, and $X \times C \subset \pi_p^{-1}(J_o)$
modulo $\eta$-null sets, 
\begin{eqnarray*}
\nu(A_{x_o}^{-1}C) 
&\geq & 
\nu((A_o \cup \{s\})^{-1}C) = \nu(A_o^{-1}C) + \nu(C) - \nu(A_o^{-1}C \cap s^{-1}C) \\
&\geq & 
m_M(I_o^{-1}J_o) + \nu(C) - \eta( A_o^{-1}(X \times C) \cap s^{-1}(X \times C)) \\
&\geq &
m_M(I_o^{-1}J_o) + \nu(C) - m_M(\tau_o(A_o)^{-1}J_o \cap \tau_o(s)^{-1}J_o) \\
&\geq &
m_M(I_o^{-1}J_o) + \nu(C) - m_M( t_o I_o^{-1}J_o \cap \tau_o(s)^{-1}J_o),
\end{eqnarray*}
or equivalently,
\[
\nu(A_{x_o}^{-1}C) - \nu(C) \geq m_M(I_o^{-1}J_o) - m_M(J_o) + m_M( \tau_o(s)^{-1} J_o \setminus t_o I_o^{-1}J_o),
\]
for all $s \in A_{x_o} \setminus A_o$.  
\end{proof}

\subsubsection{\textbf{Consequences of the overshoot bound}}

Let us retain all assumptions and notation from the previous subsection. Our aim here is to give two applications of Lemma \ref{lemma_overshoot}. \\

\textbf{Case I:} First assume that $M$ is finite,
\[
\nu(A_{x_o}^{-1}C) < \min(1,\mu(A) + \nu(C))
\]
and 
\[
m_M(I_o^{-1}J_o) = m_M(I_o) + m_M(J_o) - m_M(\{e_M\}). 
\]
By Theorem \ref{corr2} (ii) this is for instance the case if every finite quotient of $G$ is abelian. Then \eqref{overshoot} 
simplifies to 
\begin{equation}
\label{mMIo}
m_M(I_o) > \mu(A) \geq m_M(I_o) -  m_M(\{e_M\}) + m_M( \tau(s)^{-1} J_o \setminus t_o I_o^{-1}J_o),
\end{equation}
whence the rightmost term is strictly less than $m_M(\{e_M\})$ and thus zero. This readily implies the inclusion $\tau(s)^{-1}J_o \subset t_oI_o^{-1} J_o$, so if we assume that $(I_o,J_o)$ is left-balanced, then we conclude that $\tau_o(s) \in I_o t_o^{-1}$ for all $s \in A_{x_o} \setminus A_o$, and thus $A_{x_o} \subset \tau_o^{-1}(I_o t_o^{-1})$. 

Furthermore, if we denote by $G_o$ the stabilizer of $\tau^{-1}(I_o t_o^{-1})$, then this subgroup
must have finite index in $G$, and thus $m_M(\{e\}) \geq \frac{1}{[G : G_o]}$. Since $G_o A_{x_o} \subset \tau_o^{-1}(I_o t_o^{-1})$,  
we have
\[
m_M(I_o) \geq d^*(G_o A_{x_o}),
\]
and thus, by \eqref{mMIo}
\[
\mu(A) \geq m_M(I_o) - m_M(\{e_M\}) \geq d^*(G_o A_{x_o}) - \frac{1}{[G : G_o]}.
\]

\textbf{Case II:} Let us now assume that $J_o$ is equal to the closure of its interior in $M$, $(I_o,J_o)$ is left-balanced, and 
\begin{equation}
\label{eqcase}
\nu(A_{x_o}^{-1}C) = \mu(A) + \nu(C) < 1
\qand 
m_M(I_o^{-1}J_o) = m_M(I_o) + m_M(J_o).
\end{equation}
Then \eqref{overshoot} implies that 
\[
m_M(I_o) = \mu(A) \qand m_M(J_o) = \nu(C) \qand m_M( \tau(s)^{-1} J_o \setminus t_o I_o^{-1}J_o) = 0,
\]
whence $\tau_o(s)^{-1} J^o_o \subset t_o I_o^{-1} J_o$. By our assumption on $J_o$, we conclude that $\tau_o(s)^{-1}J_o \subset t_o I_o^{-1} J_o$, and since the pair is left-balanced, this implies that $\tau_o(s) \in I_o t_o^{-1}$, for all $s$ as above, whence 
\begin{equation}
\label{eqcase2}
m_M(I_o) = \mu(A) \qand m_M(J_o) = \nu(C) \qand A_{x_o} \subset \tau_o^{-1}(I_o t_o^{-1}).
\end{equation}
We are now lead to the question:
\begin{center}
\emph{When can we ensure that the conditions in Case II are satisfied?}
\end{center}
Let us recall from Proposition \ref{mainpropcont1} that if $I$ and $J$ are not sub-periodic, then
\begin{enumerate}
\item there is no finite-index subgroup $G_o < G$, extreme invariant mean $\lambda$ on $G$ and subset $A_o \subset A_{x_o}$
with $\lambda(A_o) = \mu(A)$ such that $G_o A_o \neq G$.
\item there is no finite-index subgroup $G_o < G$ such that $\nu(G_o C) < 1$.
\end{enumerate}
Let us assume that (i) and (ii) are satisfied for $A$ and $C$. By Theorem \ref{corr2} (iii), we then know that $(I,J)$ reduces to a Sturmian pair $(I_o,J_o)$ in either $M = \bT$ or $M = \bT \rtimes \{-1,1\}$. Such pairs are clearly $m_{M}$-Jordan measurable, left-balanced and $J_o$ equals the closure of its interior, so the conditions of Case II are satisfied and we get the conclusions in \eqref{eqcase2}.

\subsubsection{\textbf{Our findings}}

We shall now summarize, and slightly expand upon, our findings above in two propositions. In both of these, $G$ is assumed to be amenable, 
and our key players are:
\begin{enumerate}
\item a pointed $G$-space $(X,x_o)$ and an ergodic $\mu \in \cP_G(X)$.
\item a non-empty open $\mu$-Jordan measurable set $A \subset X$.
\item an ergodic Borel $G$-space $(Y,\nu)$ and a Borel set $C \subset Y$ with positive $\nu$-measure.
\end{enumerate}

The proof of the following proposition follows from Proposition \ref{mainpropcont1} and Case I above. 

\begin{proposition}
\label{critical}
Suppose that $\nu(A_{x_o}^{-1}C) < \min(1,\mu(A) + \nu(C))$. Then there is a subset $A_o \subset A_{x_o}$, 
a finite-index subgroup $G_o < G$ such that $G_o A_o \neq G$ and $\lambda(A_o) = \mu(A)$ for some extreme 
invariant mean $\lambda$ on $G$.

If every finite quotient of $G$ is \textsc{abelian}, then there is a finite-index subgroup $G_o < G$ such that
\[
d^*(G_o A_{x_o}) < \mu(A) + \frac{1}{[G : G_o]}.
\]
\end{proposition}

The proof of our next observation is contained under Case II above, modulo the last part, which will be proved below.

\begin{proposition}
\label{subcritical}
Suppose that $\nu(A_{x_o}^{-1}C) = \mu(A) + \nu(C) < 1$, and 
\begin{enumerate}
\item for every extremal mean $\lambda$ on $G$, for all $A_o \subset A_{x_o}$ with $\lambda(A_o) = \lambda(A_{x_o})$
and for every finite-index subgroup $G_o < G$, we have $G_o A_o = G$,

\item $\nu(G_o C) = 1$ for every finite-index subgroup $G_o < G$.
\end{enumerate}
Then $A_{x_o}$ is contained in a Sturmian set $S$ with $d^*(S) = d_*(S) = \mu(A)$. \\

Furthermore, if $G$ is 
abelian, $G \acts (Y,\nu)$ is totally ergodic and $(\bT,\tau)$ denotes the torus compactification of $G$ from 
which the Sturmian set $S$ comes, then there is a $G$-factor map $\sigma : (Y,\nu) \ra (\bT,m_{\bT})$ and a closed interval 
$J_o \subset \bT$ such that $C = \sigma^{-1}(J_o)$ modulo $\nu$-null sets.
\end{proposition}

\begin{remark}
It follows from Corollary \ref{cor_bohrset} that all invariant means assign the same value to any given 
Sturmian set $S$ in $G$, so in particular, the identity $d^*(S) = d_*(S)$ is automatic. Concerning the
last assertion in Proposition \ref{subcritical}: The assumption that $G$ is abelian is strictly speaking 
not necessary, but it simplifies the proofs significantly, which is why we choose to use it here. The general
statement (without assuming that $G$ is abelian) would also be somewhat technical to write down. 
Furthermore, total ergodicity is not necessary to assume, condition (ii) in Proposition \ref{subcritical} suffices.
\end{remark}

\subsubsection{\textbf{The proof of the last part of Proposition \ref{subcritical}}}

Let us assume that $G$ is abelian, $A_{x_o} \subset G$ satisfies (i) in Proposition \ref{subcritical} and 
$G \acts (Y,\nu)$ is totally ergodic, so that (ii) in Proposition \ref{subcritical} is automatically satisfied. 
Since $G$ is abelian, the group compactification $(K,\tau)$ is abelian as well, and thus all subgroups of $K$ 
are normal, so we may without loss of generality assume that $L$ is trivial, and thus, there is a $G$-factor map 
$\pi : (X \times Y,\eta) \ra (K,m_K)$. By Proposition \ref{mainpropcont1}, our assumptions on $A_{x_o}$
and on the action $G \acts (Y,\nu)$ force $I$ and $J$ to not be sub-periodic in $K$, and thus we are in the setting of
the conclusion \eqref{eqcase2} with $M = \bT$. In particular, if we denote by $p : K \ra \bT$ the quotient
map and $\pi_p = p \circ \pi$, then 
\begin{equation}
\label{eqcase3}
X \times C = \pi_p^{-1}(J_o), \quad \textrm{modulo $\eta$-null sets},
\end{equation}
where $J_o \subset \bT$ is a closed (proper) interval. \\

We wish to construct a $G$-factor map 
$\sigma : (Y,\nu) \ra (\bT,m_{\bT})$ such that $C = \sigma^{-1}(J_o)$
modulo $\nu$-null sets; but the way to proceed is perhaps not that clear; since $\eta$ is typically very far from 
a product measure, there is absolutely no reason to expect the $G$-factor map $\pi_p$ to "split" 
naturally into components which only depend on $x$ and $y$ respectively. For our construction, we shall 
instead use the "tightness" of \eqref{eqcase3}, combined with the fact that closed proper intervals 
$J_o$ in $\bT$ have trivial stabilizers. It will be convenient to take a bird's eye view on the matters first.

Let $(W,\theta)$ be a Borel $G$-space, and let $D \subset W$ be a Borel set. We define the map 
$\sigma_D : W \ra 2^G$ by $w \mapsto D_w = \big\{ g \in G \, : \, gw \in D \big\}$. It is not hard to 
check that this map is Borel and $G$-equivariant (recall that our action on $2^G$ is defined by 
$g.B = Bg^{-1}$). Let $U \subset 2^G$ be the clopen set $\{ B \, : \, e_G \in B \}$, and note that
$\sigma_D^{-1}(U) = D$. 
The sought-after $G$-factor map $\sigma : (Y,\nu) \ra (\bT,m_{\bT})$ above will be constructed from 
$\pi_p$ and $\sigma$'s from different Borel $G$-spaces and Borel sets therein. \\

Before we get into the construction, let us consider a special case first. If $(K,\tau)$ is a group compactification
of $G$ and $J \subset K$ is a Borel set with trivial stabilizer in $K$, that is to say, $J \neq J k $ for all $k \neq e_K$, then it
readily follows that the Borel map $\sigma_J : K \ra 2^G$ is injective, and thus has a Borel measurable inverse 
$\sigma_J^{-1}$ defined on the image of $\sigma_J$, which is also Borel (see Theorem A.4 in \cite{zi84}). In particular, 
this applies to the torus compactification $(\bT,\tau_o)$ above and the closed proper interval $J_o \subset \bT$.

Let us now turn to the construction of $\sigma : (Y,\nu) \ra (\bT,m_{\bT})$. From $(X \times Y,\eta)$ we have two natural $G$-equivariant Borel maps,
namely $\pi_p$ into $\bT$ and $\sigma_{X \times  C}$ into $2^G$. We note that 
\[
\sigma_{X \times C}(x,y) = (X \times C)_{(x,y)} = C_y = \sigma_C(y).
\]
From \eqref{eqcase3}, we have $\sigma_{(X \times C)}(x,y) = \pi_p^{-1}(J_o)_{(x,y)} = \sigma_{J_o}(\pi_p(x,y))$ for $\eta$-a.e. $(x,y)$, so
since $\sigma_{J_o}$ has a Borel measurable inverse, we can define $\sigma : (Y,\nu) \ra (\bT,m_{\bT})$ by
\begin{equation}
\label{defsigma}
\sigma(y) := (\sigma_{J_o}^{-1} \circ \sigma_C)(y) = \pi_p(x,y), \quad \textrm{for $\eta$-almost every $(x,y)$},
\end{equation}
and since $\sigma_C^{-1}(U) = C$, we note that
\[
\sigma^{-1}(J_o) = \sigma_C^{-1}(\sigma_{J_o}(J_o)) \subseteq \sigma_C^{-1}(U \cap \im \,\sigma_{J_o}) \subseteq C.
\]
Since $\sigma : Y \ra \bT$ is Borel and $G$-equivariant, and the $G$-action on $\bT$ is uniquely ergodic, the push-forward of $\nu$
must equal $m_{\bT}$. Since $m_{\bT}(J_o) = \nu(C)$, we conclude that $\sigma^{-1}(J_o) = C$ modulo 
$\nu$-null sets. 

Of course, since we also know that 
\[
A \times Y = \pi_p^{-1}(I_o), \quad \textrm{modulo $\eta$-null sets},
\]
we could have done the same thing for the set $A \subset X$, and produced a $G$-factor map 
\[
\sigma' : (X,\mu) \ra (\bT,m_{\bT}), 
\]
which again would have to coincide with $\pi_p(x,y)$ for $\eta$-almost every $(x,y)$, and $\sigma'^{-1}(I_o) = A$ modulo $\mu$-null sets. In particular, $\sigma = \sigma'$ $\eta$-almost everywhere.

The question of
which $\eta$ on $X \times Y$ that gives rise to such curious-looking identities as \eqref{defsigma} naturally arises. 
The reader is invited to check that a relatively independent joining (as well as any ergodic component thereof) over 
a common factor Borel $G$-space, which has a further factor Borel $G$-space of the form $(\bT,m_{\bT})$ will do the
job (see Section 6 in \cite{gl2003} for more details). \\

We can summarize these observations as follows. 

\begin{proposition}
\label{abeta}
Suppose that $G$ is abelian, $G \acts (Y,\nu)$ is totally ergodic, 
\[
\nu(A_{x_o}^{-1}C) = \mu(A) + \nu(C) < 1,
\]
and condition (i) in Proposition \ref{subcritical} is satisfied. Then, for every ergodic joining $\eta$ of 
$(X,\mu)$ and $(Y,\nu)$, there are 
\begin{enumerate}
\item a torus compactification $(\bT,\tau_o)$ and closed intervals $I_o, J_o \subset \bT$.
\item $G$-factor maps $\alpha : (X,\mu) \ra (\bT,m_{\bT})$ and $\beta : (Y,\nu) \ra (\bT,m_{\bT})$,
\end{enumerate}
such that $\eta\big(\big\{ (x,y) \, : \, \alpha(x) = \beta(y) \big\}) = 1$, and $A = \alpha^{-1}(I_o)$ and $C = \beta^{-1}(J_o)$
modulo null sets, where $G$ acts on $\bT$ via $\tau_o$ as in \eqref{defKtau}.
\end{proposition}

\section{Proofs of the main theorems}
\label{sec:mainresults}

\begin{adjustwidth*}{1in}{1in}
\emph{Our main \emph{dynamical} results are now rather straightforward consequences of 
Proposition \ref{critical} and Proposition \ref{subcritical}. However, in order to prove our 
main \emph{density} results, we must deal with \emph{non-ergodic} measures, which 
makes it necessary to investigate how well our ergodic-theoretical conclusions behave 
when passing to ergodic components.}
\end{adjustwidth*}

\vspace{0.2cm}

\subsection{General framework}

In what follows, let 
\begin{itemize}
\item $G$ be a countable amenable group,
\item $\xi$ an extreme invariant mean on $G$, and
\item $A$ is a subset of $G$ with $\xi(A) > 0$.
\end{itemize}
As described in Subsection \ref{subsec:hull}, we can associate
to the set $A$ a pointed $G$-space $(X,x_o)$ and (abusing notation) a \emph{clopen} subset $A \subset X$
such that $A = A_{x_o}$. From now on, we shall only work with the triple $(X,x_o,A)$, and write $A_{x_o}$
for the set in $G$. 

By Proposition \ref{onto}, the $G$-invariant probability measure $\mu := S_{x_o}^*\xi$ on $X$
is ergodic, and since $A \subset X$ is clopen, we have $\xi(A_{x_o}) = \mu(A)$ by Lemma \ref{somecons}.
\vspace{0.2cm}
\begin{center}
\emph{The notation $X,x_o,A, \xi$ and $\mu$ will be fixed throughout the rest of the section.}
\end{center}

\subsection{Proofs of the dynamical results}
\label{subsec:maindyn}
To formulate our dynamical results, we need in addition to the objects introduced above, an \emph{ergodic} 
Borel $G$-space $(Y,\nu)$ and a Borel set $B \subset Y$ with positive $\nu$-measure. The following theorem
generalizes Theorem \ref{erg1}.

\begin{theorem}
\label{ergmain1}
Suppose that
\[
\nu(A_{x_o}B) < \min\big(1,\mu(A) + \nu(B)\big).
\]
Then there exist a finite-index subgroup $G_o < G$, an extreme invariant mean $\lambda$ on $G$ and 
$A_o \subset A_{x_o}$ with $\lambda(A_o) = \mu(A)$ such that $G_o A_o \neq G$. 

If one assumes that all finite quotients of $G$ are abelian, then there exists a proper finite index subgroup $G_o < G$
such that $G_o A \neq G$, and 
\[
d^*(G_o A) < \mu(A) + \frac{1}{[G : G_o]}. 
\]
\end{theorem}

\subsubsection{\textbf{Proof of Theorem \ref{erg1} assuming Theorem \ref{ergmain1}}}

Let us assume that
\[
\nu(A_{x_o}B) < \min\big(1,\mu(A) + \nu(B)\big),
\]
and choose $\xi$ such that $\xi(A_{x_o}) = d^*(A_{x_o})$; this can be done by Proposition \ref{densextmean}. Theorem 
\ref{ergmain1} now entails that there is a finite-index subgroup $G_o$, an extreme invariant mean $\lambda$ on $G$, and
a subset $A_o \subset A_{x_o}$ with $\lambda(A_o) = \mu(A)$ such that $G_o A_o \neq G$. Since
\[
d^*(A_{x_o}) \geq d^*(A_o) \geq \lambda(A_o) = \mu(A) = \xi(A_{x_o}) = d^*(A_{x_o}),
\]
and thus $d^*(A_{o}) = d^*(A_{x_o})$, we see that $A_{x_o}$ is \emph{not} spread-out, contradicting Assumption (i) in Theorem \ref{erg1}. Assuming that all finite quotients of $G$ are abelian, Assumption (ii) in Theorem \ref{erg1} is violated along the same lines
using Theorem \ref{ergmain1}. 

\subsubsection{\textbf{Proof of Theorem \ref{ergmain1}}}

Set $C = (A_{x_o}B)^c$, and note that $A_{x_o}^{-1}C \subset B^c$. Since $1 - \nu(C) < \mu(A) + \nu(B)$ and $\nu(B) > 0$, 
we get
\[
\nu(A_{x_o}^{-1}C) \leq 1 - \nu(B) < \min(1,\mu(A) + \nu(C)),
\]
which  places us in the setting of Proposition \ref{critical}, and all properties of $A_{x_o} \subset G$ in Theorem \ref{ergmain1} readily follow from this.

\subsubsection{\textbf{Proof of Theorem \ref{erg2}}}

The following theorem generalizes Theorem \ref{erg2} in the same way that Theorem \ref{ergmain1} generalizes 
Theorem \ref{erg1}, so we omit the proof of Theorem \ref{erg2}.

\begin{theorem}
\label{ergmain2}
Suppose that
\[
\nu(A_{x_o}B) = \mu(A) + \nu(B) < 1.
\]
If $A_{x_o}$ is spread-out, and $A_{x_o}B$ does not contain, modulo $\nu$-null sets, a Borel set with positive measure which is invariant under
a finite-index subgroup $G_o < G$, then $A_{x_o}$ is contained in a Sturmian set $S$ with $\xi(S) = \xi(A_{x_o}) = \mu(A)$.\\

Furthermore, if one in addition assumes that $G$ is abelian and $G \acts (Y,\nu)$ is totally ergodic, then, for every ergodic 
joining $\eta$ of $(X,\mu)$ and $(Y,\nu)$, there exist
\begin{enumerate}
\item a torus compactification $(\bT,\tau_o)$ and closed intervals $I_o, J_o \subset \bT$.
\item $G$-factor maps $\alpha : (X,\mu) \ra (\bT,m_{\bT})$ and $\beta : (Y,\nu) \ra (\bT,m_{\bT})$,
\end{enumerate}
such that 
\begin{itemize}
\item $\eta(\{(x,y) \, : \, \alpha(x) = \beta(y) \big\}) = 1$, 
\item $A = \alpha^{-1}(I_o)$ and $B = \beta^{-1}(J_o)$ modulo null sets,
\end{itemize}
where $G$ acts on $\bT$ via $\tau_o$ as in \eqref{defKtau}.
\end{theorem}

\subsubsection{\textbf{Proof of Theorem \ref{ergmain2}}}

Set $C = (A_{x_o}B)^c$, and note that $A_{x_o}^{-1}C \subset B^c$. Since $1 - \nu(C) = \mu(A) + \nu(B) < 1$ and $\nu(B) > 0$, 
we get $\nu(A_{x_o}^{-1}C) \leq 1 - \nu(B) = \mu(A) + \nu(C) < 1$. Since $A_{x_o}$ is assumed to be spread-out, we see by
Proposition \ref{critical} that the first inequality cannot be strict, whence
\[
\nu(A_{x_o}^{-1}C) = \mu(A) + \nu(C) < 1.
\]
which  places us in the setting of Proposition \ref{subcritical} and Proposition \ref{abeta}. \\

Assume that there exists a finite-index subgroup $G_o < G$ such
that $\nu(G_o C) < 1$. Then, by taking complements, we see that the Borel set $D := \bigcap_{g \in G_o} gA_{x_o} B \subset Y$ has positive $\nu$-measure and is $G_o$-invariant, contradicting our second assumption in Theorem \ref{ergmain2}. Hence the 
conditions of Proposition \ref{subcritical} are satisfied, and the conclusions about the set $A_{x_o}$ follow. \\

Let us now assume that $G$ is abelian and $G \acts (Y,\nu)$ is totally ergodic. Let us also fix an ergodic joining $\eta$ of $(X,\mu)$
and $(Y,\nu)$. Since $\nu(A_{x_o}^{-1}C) = \mu(A) + \nu(C) < 1$, Proposition \ref{abeta} ensures the existence of a torus compactification $(\bT,\tau_o)$ and $G$-factor maps 
\[
\alpha : (X,\mu) \ra (\bT,m_{\bT}) \qand \beta : (Y,\nu) \ra (\bT,m_{\bT})
\]
such that $\eta(\{(x,y) \, : \, \alpha(x) = \beta(y) \big\}) = 1$ holds, and closed intervals $I_o, H_o \subset \bT$ such that
\[
A = \alpha^{-1}(I_o) \qand C = \beta^{-1}(H_o),
\]
modulo null sets. Furthermore, upon going into the arguments of Proposition \ref{subcritical}, we see that there is an 
element $t_o \in \bT$ such that $A_{x_o} \subset \tau_o^{-1}(I_o t_o^{-1})$. It follows from chain of identities above that 
\[
A_{x_o}^{-1}C = B^c \subset \tau_o^{-1}(t_o I^{-1}) \beta^{-1}(H_o),
\]
modulo null sets. It is also not hard to see that since both $I_o$ and $H_o$ are intervals and $\beta$ is $G$-equivariant, we have 
\[
\tau_o^{-1}(t_o I^{-1}) \beta^{-1}(H_o) = \beta^{-1}(t_o I_o^{-1}H_o)
\]
modulo null sets. Since $J_o := (t_o I_o^{-1} H_o)^c$ is again an interval in $\bT$ with
\[
m_{\bT}(J_o) = 1 - m_{\bT}(I_o^{-1}H_o) = 1 - m_{\bT}(I_o) - m_{\bT}(H_o) = \nu(B),
\]
we conclude that $B = \beta^{-1}(J_o)$ modulo null sets, which finishes the proof.

\begin{remark}
Note that in both proofs above, we went from action sets of the form $A_{x_o} B$ to action sets of the form $A_{x_o}^{-1}C$. 
This was of course done so that the results in the previous section could be applied, but it is natural to ask whether this 
conversion is necessary - surely this is not the case for \emph{abelian} $G$. The need for this twist can be traced 
to Lemma \ref{symmetry}; in the proof of this lemma,we heavily use that the map $x \mapsto \nu(A_{x}^{-1}C)$ is $G$-invariant. 
This is not the case for the map $x \mapsto \nu(A_x B)$, unless of course $G$ is abelian.  
\end{remark}

\subsection{Proofs of the density results}
\label{subsec:maindens}

We retain the notation for $X,x_o,A,\xi$ and $\mu$ introduced in the beginning of the section. We shall further 
fix a subset $B \subset G$, and add conditions on it as we go along. We associate to $B$ a pointed $G$-space 
$(Y,y_o)$ so that there is a clopen set $B \subset Y$ (abuse of notation) such that $B = B_{y_o}$, where the set $B$ on 
the left is in $G$, and the set $B$ on the right is the clopen set in $Y$. To avoid this abuse of notation,
we shall from now on only refer to the set in $G$ as $B_{y_o}$.

\subsubsection{\textbf{Proofs of Theorem \ref{dens1ubd} and Theorem \ref{dens2ubd}}}

By Proposition \ref{densextmean}, we can find \emph{extreme} invariant means $\lambda_{+}$ and $\lambda_{-}$
on $G$ such that
\begin{equation}
\label{oeta}
\lambda_{+}(B_{y_o}) = d^*(B_{y_o})
\qand 
\lambda_{+}(A_{x_o}B_{y_o}) \leq d^*(A_{x_o} B_{y_o})
\end{equation}
and
\begin{equation}
\label{ueta}
\lambda_{-}(B_{y_o}) \geq d_{*}(B_{y_o})
\qand 
\lambda_{-}(A_{x_o}B_{y_o}) = d_*(A_{x_o} B_{y_o}).
\end{equation}
By Proposition \ref{onto}, the $G$-invariant probability measures $\nu_{+} = S_{y_o}^*\lambda_{+}$
and $\nu_{-} = S_{y_o}^*\lambda_{-}$ are \emph{ergodic}, and by Lemma \ref{somecons} we have
\begin{equation}
\label{oeta}
\nu_{+}(B) = d^*(B_{y_o})
\qand 
\nu_{+}(A_{x_o}B) \leq d^*(A_{x_o} B_{y_o})
\end{equation}
and
\begin{equation}
\label{ueta}
\nu_{-}(B) \geq d_{*}(B_{y_o})
\qand 
\nu_{-}(A_{x_o}B_{y_o}) \leq d_*(A_{x_o} B_{y_o}).
\end{equation}
If we enforce the assumptions in Theorem \ref{dens1ubd} on the set $A_{x_o}$, then Theorem \ref{erg1} readily implies
Theorem \ref{dens1ubd}. \\

Towards the proof of Theorem \ref{dens2ubd}, suppose that $A_{x_o}$ is spread-out, $B_{y_o}$ large and $A_{x_o} B_{y_o}$ 
does not contain a piecewise periodic set. If $d^*(A_{x_o} B_{y_o}) = d^*(A_{x_o}) + d^*(B_{y_o}) < 1$, then, with the notation 
above,
\[
\nu_{+}(A_{x_o}B) \leq d^*(A_{x_o}) + \nu_{+}(B) 
\]
Since $A_{x_o}$ is spread-out and $\nu_{+}$ is ergodic, Theorem \ref{erg1} shows that the inequality cannot be strict and
thus $\nu_{+}(A_{x_o}B) = d^*(A_{x_o}) + \nu_{+}(B) < 1$. At this point, Theorem \ref{erg2} tells us $A_{x_o}$ is contained in a 
Sturmian set with the same upper Banach density as $A_{x_o}$, \emph{provided} that $A_{x_o}B$ does not contain a Borel 
set $Z$ which is invariant under a finite-index subgroup $G_o$. To prove that this is not the case, we argue by contradiction, 
and apply Lemma \ref{convtolemma} to the ergodic $G$-space $(Y,\nu)$ with $Q = Z$ and $U = A_{x_o} B \subset Y$. We conclude that there is a non-empty $G_o$-invariant set $Q_o \subset G$ and a thick set $T \subset G$ such that
\[
U_{y_o} = A_{x_o} B_{y_o} \supset Q_o \cap T,
\]
contradicting our assumption that $A_{x_o} B_{y_o}$ does not contain a piecewise periodic set.

\subsubsection{\textbf{How to deal with asymptotic densities}}
\label{subsec:dealwithit}
As we have just seen, Theorem \ref{dens1ubd} and Theorem \ref{dens2ubd} are rather direct consequences of
Theorem \ref{erg1} and Theorem \ref{erg2}. The main reason for this is that the measures $\nu_{+}$ and $\nu_{-}$ on $Y$
that we end up with are \emph{ergodic} - or, equivalently, the maximixing/minimizing invariant means $\lambda_{+}$ and $\lambda_{-}$
are extreme points in $\cL_G$. This will no longer be the case when we study asymptotic densities. \\

Towards the proofs of Theorem 
\ref{dens1} and Theorem \ref{dens2}, we shall begin by fixing a F\o lner sequence $(F_n)$ in $G$
once and for all. We can then use \eqref{asympdensmean} and \eqref{besureof} to produce (not necessarily 
extreme) invariant means $\overline{\lambda}$ and $\underline{\lambda}$ such that 
\[
\overline{\lambda}(B) = \overline{d}_{(F_n)}(B_{y_o})
\qand 
\overline{\lambda}(A_{x_o}B_{y_o}) \leq \overline{d}_{(F_n)}(A_{x_o} B_{y_o})
\]
and
\[
\underline{\lambda}(B) \geq \underline{d}_{(F_n)}(B_{y_o})
\qand 
\underline{\lambda}(A_{x_o}B_{y_o}) = \underline{d}_{(F_n)}(A_{x_o} B_{y_o}).
\]
If we write $\overline{\nu} = S_{y_o}^*\overline{\lambda}$ and $\underline{\nu} = S_{y_o}^*\underline{\lambda}$, then
by Lemma \ref{somecons}, 
\begin{equation}
\label{onu}
\overline{\nu}(B) =  \overline{d}_{(F_n)}(B_{y_o})
\qand
\overline{\nu}(A_{x_o}B) \leq \overline{d}_{(F_n)}(A_{x_o} B_{y_o}),
\end{equation}
and
\begin{equation}
\label{unu}
\underline{\nu}(B) \geq  \underline{d}_{(F_n)}(B_{y_o})
\qand
\underline{\nu}(A_{x_o}B) \leq \underline{d}_{(F_n)}(A_{x_o} B_{y_o}).
\end{equation}
It will be useful in the subsequent arguments to recast our assumptions on $B_{y_o}$ and $AB_{y_o}$ in Theorem \ref{dens1} and Theorem
\ref{dens2} as properties of sets in $Y$. The following lemma does so. 

\begin{lemma}
\label{conversion}
With the notation and assumptions above,
\begin{enumerate}
\item $B_{y_o}$ is syndetic $\iff$ $\nu(B) > 0$ for all $\nu \in \cP_G(Y)$.
\item $A_{x_o}B_{y_o}$ is not thick $\implies$ $\nu(A_{x_o}B) < 1$ for all $\nu \in \cP_G(Y)$.
\item $A_{x_o} B_{y_o}$ does not contain a piecewise periodic set $\implies$ For every $\nu \in \cP_G^{\textrm{erg}}(Y)$, $A_{x_o}B$ does not 
contain a Borel set with positive $\nu$-measure which is invariant under a finite-index subgroup.
\end{enumerate}
\end{lemma}

\begin{proof}
(i) and (ii) are immediate consequences of Lemma \ref{somecons} and Lemma \ref{thicksyn}. 

(iii) Assume that there exists $\nu \in \cP_G^{\textrm{erg}}(Y)$ and a Borel set $Z \subset A_{x_o} B$ with positive $\nu$-measure, 
which is invariant under a finite-index subgroup. Apply Lemma \ref{convtolemma} to the ergodic $G$-space $(Y,\nu)$ with $Q = Z$ and $U = A_{x_o} B$. We conclude that there is a $G_o$-invariant set $Q_o \subset G$ and a thick set $T \subset G$ such that
\[
U_{y_o} = A_{x_o} B_{y_o} \supset Q_o \cap T,
\]
showing in particular that $A_{x_o} B_{y_o}$ contains a piecewise periodic set. 
\end{proof}

\subsubsection{\textbf{Ergodic decompositions}}

In what follows, let $\nu$ be a $G$-invariant (not necessarily ergodic) Borel probability measure on $Y$.
In the applications that will follow, we will consider the cases $\nu = \overline{\nu}$ and $\nu = \underline{\nu}$. It is a 
standard fact in ergodic theory (see for instance Theorem 4.8 in \cite{eiwa2015}), that one can decompose
$\nu$ into \emph{ergodic components}, that is to say, there exists a probability measure $\kappa$ on 
$\cP_G(Y)$, which is concentrated on the set of ergodic measures, such that
\begin{equation}
\label{ergcompo}
\nu(D) = \int_{\cP^{\textrm{erg}}(Y)} \nu'(D) \, d\kappa(\nu'), \quad \textrm{for all Borel sets $D \subset Y$}.
\end{equation}

\subsubsection{\textbf{Proof of Theorem \ref{dens1}}}
\label{prfdens1}

Recall that $\nu$ is fixed. Let us assume that 
\begin{enumerate}
\item $A_{x_o}$ is spread-out (or, if every finite quotient of $G$ is abelian, that \eqref{perbnd} does
not hold for any finite-index subgroup $G_o$).
\item $B_{y_o}$ is syndetic; Lemma \ref{conversion}  implies that
$\nu'(B) > 0$ for all $\nu' \in \supp(\kappa)$.
\item $A_{x_o} B_{y_o}$ is not thick; Lemma \ref{conversion}  implies that
$\nu'(A_{x_o}B) < 1$ for all $\nu' \in \supp(\kappa)$.
\end{enumerate}
By Theorem \ref{erg1}, applied to each ergodic component $\nu'$ of $\nu$, we can now conclude
\[
1 > \nu'(A_{x_o}B) \geq \mu(A) + \nu'(B),
\]
and thus, by \eqref{ergcompo},
\begin{eqnarray*}
\nu(A_{x_o}B) 
&=& 
\int_{\cP^{\textrm{erg}}(Y)} \nu'(A_{x_o}B) \, d\kappa(\nu') \\
&\geq &
\int_{\cP^{\textrm{erg}}(Y)}  (\mu(A) + \nu'(B)) \, d\kappa(\nu') \\
&=&
\mu(A) + \nu(B).
\end{eqnarray*}
Let us now pick $\mu \in \cP^{\textrm{erg}}_G(X)$ so that $d^*(A_{x_o}) = \mu(A)$. Theorem \ref{dens1} readily follows if we apply
the previous inequalities to $\nu = \overline{\nu}$ and $\nu = \underline{\nu}$ respectively, together with \eqref{onu} and \eqref{unu}.

\subsubsection{\textbf{Proof of Theorem \ref{dens2}}}
\label{prfdens2}

Let us now assume that 
\begin{enumerate}
\item $A_{x_o}$ is spread-out.
\item $B_{y_o}$ is syndetic; Lemma \ref{conversion}  implies that
$\nu'(B) > 0$ for all $\nu' \in \supp(\kappa)$.
\item $A_{x_o} B_{y_o}$ does not contain a piecewise periodic subset; Lemma \ref{conversion} implies that
for every $\nu' \in \supp(\kappa)$, the action set $A_{x_o}B$ does not contain a Borel set with positive $\nu'$-measure which 
is invariant under a finite-index subgroup.
\item $\nu(A_{x_o} B) = \mu(A) + \nu(B) < 1$.
\end{enumerate}

We claim that
\[
\nu'(A_{x_o}B) = \mu(A) + \nu'(B) < 1, \quad \textrm{for $\kappa$-a.e. $\nu'$}.
\]
Indeed, by Theorem \ref{erg1} applied to each $\nu'$ (using the assumption that $A_{x_o}$ is spread-out), we know
that $\nu'(A_{x_o}B) \geq \mu(A) + \nu'(B)$ for $\kappa$-a.e. $\nu'$. However, by \eqref{ergcompo} and our assumption (iv)
above, we also have
\[
\int_{\cP(Y)} \nu'(A_{x_o} B) \, d\kappa(\nu') = \nu(A_{x_o} B) = \mu(A) + \nu(B) = \int_{\cP(Y)} \big(\mu(A) + \nu'(B) \big) d\kappa(\nu'),
\]
so the inequality $\nu'(A_{x_o}B) \geq \mu(A) + \nu'(B)$ cannot be strict on a set of positive $\kappa$-measure.
%

We now see that
all conditions of Theorem \ref{erg2} are satisfied for every ergodic component $\nu'$, and
thus we conclude that $A_{x_o}$ is contained in a Sturmian set $S$ with upper Banach density
equal to $\mu(A)$. When $\mu$ (or equivalently, $\xi$) is chosen so that $\xi(A_{x_o}) = d^*(A_{x_o})$, 
Theorem \ref{dens2} follows.

\addtocontents{toc}{\protect\setcounter{tocdepth}{1}}
\subsection{Some auxiliary consequences of our arguments (optional)}
\addtocontents{toc}{\protect\setcounter{tocdepth}{2}}
\label{subsec:aux}

Let us now present the promised corollaries of Theorem \ref{corr1}.

\begin{theorem}
\label{addcons1}
Let $G$ be a finitely generated amenable group, and suppose that $G \acts (Y,\nu)$ is an ergodic Borel $G$-space. 
Let $A \subset G$ be a large set and $B \subset Y$ a Borel set with positive $\nu$-measure. 
\begin{enumerate}
\item If $G$ is simple, then $\nu(AB) = 1$.
\item If $G$ is a torsion group and $A \subset G$ is spread-out, then $\nu(AB) = 1$.
\end{enumerate}
\end{theorem}

\begin{proof}
Let $(X,x_o)$ be the hull associated to $A$, abuse notation (as many times before) and denote by $A$ the clopen set in $X$ 
so that $A = A_{x_o}$. Fix a Borel set $B \subset Y$ with positive $\nu$-measure, and set $C = (A_{x_o}B)^c$. If we assume that 
$\nu(A_{x_o}B) < 1$, then
\[
\nu(C) > 0 \qand \nu(A_{x_o}^{-1}C) < \nu(B^c) < 1.
\]
Let $I$ and $J$ denote the sets in Theorem \ref{corr1}. 

(i) If $G$ is simple, then Lemma \ref{compactify} (iv) shows that any group compactification $(K,\tau)$ of $G$ is trivial,
and thus $1 = m_K(I^{-1}J) \leq \nu(A_{x_o}^{-1}C)$ by Theorem \ref{corr1}. 

(ii) If $G$ is a torsion group, then Lemma \ref{compactify} (iii) shows that any group compactification $(K,\tau)$ of $G$ is totally 
disconnected. If $A_{x_o}$ is spread-out, then it follows from Proposition \ref{mainpropcont1} that the set $I \subset K$ cannot
be sub-periodic. In particular, for every open subgroup $U < K$, we have $UI = K$. Since the open subgroups form a neighborhood
basis of the identity, this implies in turn that $I$ must be dense in $K$, and thus $m_K(I^{-1}J) = 1$ by Theorem \ref{potpurri} (i),
which, via Theorem \ref{corr1}, shows that $\nu(A_{x_o}^{-1}C) = 1$ 
\end{proof}

\begin{remark}
In fact, the proof of (ii) gives a bit more: If $G$ is a finitely generated amenable torsion group and $G \acts (Y,\nu)$ is
a (non-trivial) totally ergodic Borel $G$-space, then it is weakly mixing. Indeed, if $G \acts (Y,\nu)$ is \emph{not} weakly mixing, 
then by Lemma \ref{Mackey}, there exist 
\begin{enumerate}
\item a metrizable group compactification $(K,\tau)$ of $G$, and a closed \emph{proper} subgroup $L < K$.
\item a $G$-factor map $\pi : (Y,\nu) \ra (K/L,m_{K/L})$.
\end{enumerate} 
By assumption, every finite-index subgroup $G_o < G$ acts ergodically on $(Y,\nu)$ and thus also on $(K/L,m_{K/L})$. 
It is not hard to see that this implies that $UL = K$ for every \emph{open} subgroup $U < K$. By Lemma \ref{compactify} (iii),
$K$ is totally disconnected, so we can find a decreasing chain $(U_n)$ of open subgroups such that $\bigcap_n U_n = \{e_K\}$.
Since $L$ is proper, there exists $t \in K$ such that $tL \cap L = \emptyset$, and since $\bigcap_n tL \cap U_n = \emptyset$, we
can by compactness find $n$ such that $tL \cap U_n = \emptyset$, and thus $t \notin U_n L$, which is a contradiction. We 
conclude that $L = K$, and thus $K/L$ is trivial, whence $G \acts (Y,\nu)$ is weakly mixing.
\end{remark}

The following corollary is immediate.

\begin{corollary}
\label{addcons2}
Let $G$ be a finitely generated amenable group, and let $A \subset G$ be a large set. Then $d^*(AB) = 1$ for every 
large set $B \subset G$ if either $G$ is simple, or if $G$ is a torsion group and $A$ is spread-out.
\end{corollary}

\begin{remark}
The first assertion (when $G$ is simple) was essentially proved by Bergelson and Furstenberg in \cite{befu09} - they prove the same
result under the assumption that $G$ is minimally almost periodic, that is to say, $G$ admits no non-trivial group 
compactification (this is the only property that we use as well).
\end{remark}

\section{Counterexample machine for semi-direct products}
\label{sec:cntex}
\begin{adjustwidth*}{1in}{1in}
\emph{We develop a "machine" which supplies counterexamples to
certain conjectural "symmetrized" versions of our main results concerning upper 
Banach densities of product sets in groups which are far 
from being abelian.}
\end{adjustwidth*}

\vspace{0.2cm}

\subsection{General setting}

Throughout this section, let $G$ be a countable group which is a product of two distinguished 
subgroups $N$ and $L$, where $N$ is \emph{abelian} and \emph{normal} in $G$, and $N \cap L = \{e_G\}$. 
We shall assume that there is a proper \emph{finitely generated} subgroup $\Lambda$ of $N$ with 
the property that for \emph{every} finite subset $F \subset N$, there is an element $l \in L$ such that 
$l F l^{-1} \subset \Lambda$. The reader can verify that these assumptions imply that
\begin{itemize}
\item $G$ is amenable $\iff$ $L$ is amenable.
\item $N$ is \emph{not} finitely generated. 
\end{itemize}

The two main examples to keep in mind are
\begin{equation}
\label{cntex1}
G = \bZ[1/p] \rtimes \langle p \rangle \qand L = \langle p \rangle \qand N = \bZ[1/p] \qand \Lambda = \bZ,
\end{equation}
for some prime number $p$, which acts by multiplication on $\bZ[1/p]$, and
\begin{equation}
\label{cntex2}
G = \bQ \rtimes \bQ^* \qand L = \bQ^* \qand N = \bQ \qand \Lambda = \bZ.
\end{equation}
In both of these examples, the group $G$ is two-step solvable, and hence amenable. \\

We shall from now on assume that $L$, and hence $G$, is amenable. Also, to avoid confusion, we 
denote by $d^*_G$ and $d^*_L$ the upper Banach densities on $G$ and $L$ respectively. 
Towards the proofs of Theorem \ref{thm_counterex1} and Theorem \ref{thm_counterex2} we 
record in the next proposition some peculiar behaviors of the sets
\begin{equation}
\label{defST}
S = L \Lambda \qand T = (S^{-1}S)^c
\end{equation}
with respect to the upper Banach density on $G$.

\begin{proposition}
\label{counterexmachine}
With the notation and assumptions above, we have
\begin{enumerate}
\item $d_G^*(S) = d_G^*(T) = 1$.
\item For any $A_o, B_o \subset L$, we have $d_G^*(AB) \leq d^*_L(A_o B_o)$, where
\[
A = NA_o \cap S \qand B = (NB_o \cap T) \sqcup \{e_G\}.
\]
\end{enumerate}
\end{proposition}

Property (i) will be proved below. Towards the proof of (ii), we note that 
\begin{eqnarray*}
AB
&=& 
(NA_o \cap S)(NB_o \cap T) \cup (NA_o \cap S) \\
&\subseteq &
(NA_o B_o \cap ST) \sqcup (NA_o \cap S).
\end{eqnarray*}
Furthermore, since $S$ is left $L$-invariant and every element $g \in G$ can be written on the form $nl$ for some 
$n \in N$ and $l \in L$, we see that for every invariant mean $\lambda$ on $G$ and $g \in G$, 
\[
\lambda(NA_o B_o \cap g ST) = \lambda(NA_o B_o \cap ST)
\qand
\lambda(NA_o \cap g S) = \lambda(NA_o \cap S).
\]
In particular, if we assume that $\lambda$ is extreme in $\cL_G$, then Proposition \ref{wet} implies
that
\begin{equation}
\label{blaj}
\lambda(NA_o B_o \cap ST) = \lambda(NA_o B_o) \lambda(ST)
\qand
\lambda(NA_o \cap S) = \lambda(NA_o) \lambda(S).
\end{equation}
We note that $\lambda_o( \cdot ) = \lambda(N \cdot)$ defines an invariant mean on $L$ (this is simply the push-forward of
$\lambda$ under the quotient map $G \ra G/N$), and thus 
$\lambda(NA_o B_o) \leq d^*_L(A_o B_o)$ by Proposition \ref{densextmean}. Hence, for every extreme 
invariant mean $\lambda$ on $G$, 
\begin{eqnarray*}
\lambda(AB) 
&\leq & 
\lambda(NA_o B_o \cap ST) + \lambda(NA_o \cap S) \\
&=&
\lambda(NAB_o)\lambda(ST) + \lambda(NA_o) \lambda(S) \\
&\leq &
\lambda(NA_o B_o) \, \lambda(ST \sqcup S) \\
&\leq & 
\lambda(NA_o B_o) \\
&\leq &
d^*_L(A_o B_o),
\end{eqnarray*}
where the first identity follows from \eqref{blaj} and second inequality follows from monotonicity of $\lambda$ and the fact that $ST \cap S = \emptyset$. Finally, 
Proposition \ref{densextmean} allows us to pick an extreme $\lambda$ in $\cL_G$ such that $\lambda(AB) = d^*_G(AB)$, which finishes the proof of (ii). 

Let us now turn to the proof of (i). To prove that $d_G^*(S) = 1$ it suffices by Lemma \ref{thicksyn} (i) to show that
for every finite subset $F \subset G$, there exists $g \in G$ such that $Fg \subset S$. In our setting, we may without 
loss of generality assume that $F$ is of the form $F_LF_N$ where $F_L \subset L$ and $F_N \subset N$ are finite 
sets. By our assumptions on $L,N$ and $\Lambda$, we can now find $l \in L$ such that $l F_N l^{-1} \subset \Lambda$,
and thus
\[
Fl^{-1} = (F_L l^{-1})(lF_Nl^{-1}) \subset L \Lambda = S,
\]
which finishes the proof that $d_G^*(S) = 1$. Towards the proof of the second claim, we note that in order to show that $d_G^*(T) = 1$, or 
equivalently, $d^G_*(S^{-1}S) = 0$, it suffices by Lemma \ref{thicksyn} (ii) to show that there is \emph{no} 
finite subset $F \subset G$ such that $FS^{-1}S = F\Lambda L \Lambda = G$; in particular, it would be enough
to show that there is no finite set of the form $F_N F_L$, where $F_N \subset N$ and $F_L \subset L$ are 
finite subsets, such that 
\[
F_N F_L \Lambda L \Lambda \cap N = N.
\]
To reach a contradiction, let us assume 
that such sets $F_N$ and $F_L$ exist, and note that since the intersection $N \cap L$ is trivial, this implies that
\begin{equation}
\label{FNFL}
F_N F_L \Lambda L \Lambda \cap N = F_N\Big( \bigcup_{l \in F_L} l \Lambda l^{-1}\Big) \Lambda = N.
\end{equation}
Indeed, since $L \cap N = \{e\}$ and $\Lambda < N$, the only elements in the set $F_N F_L \Lambda L$ which belong to $N$ are the ones of the form
\[
f_N f_L \lambda_1 f_L^{-1} \lambda_2,
\]
where $f_N \in F_N$, $f_L \in F_L$ and $\lambda_1, \lambda_2 \in \Lambda$. Conversely, every element of this form belongs to the intersection 
$F_N F_L \Lambda L \Lambda \cap N$. 

By assumption, $\Lambda$ is generated by some \emph{finite} set $Q_o$, so we conclude from \eqref{FNFL} that $N$
must be generated by the finite set $F_N \cup \bigcup_{l \in F_L} l Q_o l^{-1} \cup Q_o$. However, we observed
already in the beginning of the section that under our assumptions, $N$ cannot be finitely generated. This finishes the 
proof of (i).

\subsection{Constructing counterexamples}

The constructions of the counterexamples in Theorem \ref{thm_counterex1} and Theorem \ref{thm_counterex2} share the same 
basic structure with another. Let $(K,\tau_o)$ be a metrizable
group compactification of $L$; since $N$ is normal, we may extend $(K,\tau_o)$ to a group compactification $(K,\tau)$ of $G$ by setting 
$\tau(nl) = \tau_o(l)$ for $n \in N$ and $l \in L$. Let $I, J \subset K$ be closed $m_K$-Jordan measurable subsets
which are equal to the closures of their interiors, and such that $IJ$ is $m_K$-Jordan measurable as well. Set
\begin{equation}
\label{basicform}
A = \tau^{-1}(I) \cap S \qand B = (\tau^{-1}(J) \cap T) \sqcup \{e_G\}.
\end{equation}
We note that $A$ and $B$ are constructed from the sets $A_o = \tau_o^{-1}(I)$ and $B_o = \tau_o^{-1}(J)$ exactly 
as in Proposition \ref{counterexmachine}, whence
\begin{equation}
\label{ineqcntex}
d^*_G(AB) \leq d^*_L(A_o B_o) \leq d^*_L(\tau_o^{-1}(IJ)) = m_K(IJ),
\end{equation} 
where we in the last equality used Corollary \ref{cor_bohrset} and our assumption $IJ$ is $m_K$-Jordan measurable. 
The lemma below records some further important properties of these sets.

\begin{lemma}
\label{auxcnt}
With the notation and assumptions above, we have:
\begin{enumerate}
\item $d^*_G(A) = m_K(I)$ and $d^*_G(B) = m_K(J)$
\item Let $G_o < G$ be a finite-index subgroup:
\begin{enumerate}
\item If $\tau(G_o)J = K$, then $G_o B = G$.
\item If $J$ has trivial stabilizer and $\tau(G_o) \cap J = \emptyset$, then
\[
d^*_G(G_o B) > d_G^*(B) + \frac{1}{[G : G_o]}.
\]
\end{enumerate}
\item If $K$ is connected, then $A$ and $B$ are spread-out in $G$, and, if in addition, $I J \cup I \neq K$,
then $AB$ does not contain a piecewise periodic set.
\item If $e_K \notin J$, then $B$ is not contained in a Sturmian set in $G$ with the same upper Banach density as $B$.
\end{enumerate}
\end{lemma}

\begin{proof}
(i) Since $d^*_G(S) = d^*_G(T) = 1$, we can by Proposition \ref{densextmean} find invariant means $\lambda_1$
and $\lambda_2$ on $G$ such that $\lambda_1(S) = \lambda_2(T) = 1$, and thus, in combination with Corollary
\ref{cor_bohrset}, 
\[
\lambda_1(A) = \lambda_1(\tau^{-1}(I)) = m_K(I)
\qand
 \lambda_2(A) = \lambda_2(\tau^{-1}(J)) = m_K(J).
\]
Clearly, $A \subset \tau^{-1}(I)$ and $B \subset \tau^{-1}(J) \cup \{e_G\}$, whence $d_G^*(A) \leq m_K(I)$ and 
$d_G^*(B) \leq m_K(J)$, which finishes the proof in view of Proposition \ref{densextmean}.

(ii) Since $T$ is thick, Lemma \ref{stuff} (ii) tells us that 
$G_o(\tau^{-1}(J) \cap T) = \tau^{-1}(\tau(G_o)J)$, and thus $G_o B = \tau^{-1}(\tau(G_o)J) \cup G_o$. This
finishes (a). If the conditions in (b) hold, then we claim
that the union is disjoint; if not, $J \cap \tau(G_o) \neq \emptyset$, contradicting our assumption. Hence, for 
any invariant mean $\lambda$ on $G$, 
\[
\lambda(G_oB) \geq \lambda(\tau^{-1}(\tau(G_o)J)) + \lambda(G_o) > \lambda(\tau^{-1}(J)) + \frac{1}{[G : G_o]} = m_K(J) + \frac{1}{[G:G_o]},
\]
where we in the last equality used Corollary \ref{cor_bohrset}, and in the second inequality the fact that $m_K$-Jordan
measurability and trivial stabilizer of $J$ implies that $J' = \tau(G_o) J \setminus J$ contains a non-empty open set in 
$K$, and thus $\lambda(\tau^{-1}(J'))$ is strictly positive by Lemma \ref{somecons} (iii). The fact that 
$\lambda(G_o) = 1/[G : G_o]$ is left to the reader.

(iii) First note that $AB \subset \tau^{-1}(IJ \cup I)$, so if $AB$ contains $Q \cap U$ for some right $G_o$-invariant set $Q$ and 
thick set $U$, then $\overline{\tau(Q \cap U)} \subset IJ \cup I \neq K$. However, this contradicts Lemma \ref{stuff} (iii).
It thus remains to show that $A$ and $B$ are spread-out. To do this, first note that if $A' \subset A$ has the same upper
Banach density, then we can write $A' = \tau^{-1}(I) \cap U_I$ for some thick set $U_I \subset G$. We wish to prove that $G_o A' = G$
for any finite-index subgroup $G_o$ of $G$.
By Lemma \ref{stuff} (ii), $G_o A' = \tau^{-1}(\tau(G_o)I) = G$, where the last identity follows since $K$ is connected and thus the 
image of $G_o$ under $\tau$ is dense in $K$. Similarly, any subset $B' \subset B$ with
the same upper Banach density can be written of the form $\tau^{-1}(J) \cap U_J$ for some thick set $U_J$, possibly 
adding $e_G$ depending on whether it belongs to $B'$ or not. Again, by Lemma \ref{stuff} (ii) and the fact that $K$ is connected, 
$G_o B' = G$, so $B$ is spread-out.

(iv) Suppose that $(M,\theta)$ is a group compactification of $G$ and $J' \subset M$ a closed $m_M$-Jordan measurable 
subset equal to the closure of its interior such that 
\[
m_M(J') = d_G^*(\theta^{-1}(J')) = d_G^*(B) = m_K(J) \qand B = (\tau^{-1}(J) \cap T) \sqcup \{e_G\} \subset \theta^{-1}(J').
\]
We wish to show that such $M, \tau$ and $J'$ cannot exist, proving in particular that $B$ cannot be contained in a Sturmian
set with the same upper Banach density. To disprove existence, assume that these things exist, and
consider the homomorphism $\xi : G \ra K \times M$ defined by $\xi(g) = (\tau(g),\theta(g))$, and denote by $E$ the closure
of $\xi(G)$ in $K \times M$. Then $(E,\xi)$ is a group compactification of $G$ and $E$ maps onto both $K$ and $M$.
 We set
\[
C = (J \times M) \cap E \qand D = (K \times J') \cap E,
\]
and leave it to the reader to show that $C$ and $D$ are closed, $m_E$-Jordan measurable and equal to the 
closures of their interiors in $E$. Furthermore, $m_E(C) = m_K(J) = m_M(J') = m_E(D)$. Since 
$\xi^{-1}(C) = \tau^{-1}(J)$ and $\xi^{-1}(D) = \theta^{-1}(J')$, we see that $\xi^{-1}(C) \sqcup \{e_G\} \subset \xi^{-1}(D)$,
whence
\[
C = \overline{C^o} \subset \overline{C \cap \xi(G)} \subset D.
\]
In particular, since $C$ and $D$ are $m_E$-Jordan measurable and $m_E(C) = m_E(D)$, the open set 
$D^o \setminus C$ is null, and thus empty, whence $D^o \subset C$. Since $D$ equals the closure of its
interior, we conclude that $C = D$. However, going back a few lines, we see that this implies that $e_E \in C$,
and thus $e_K \in J$, which contradicts our assumption.
\end{proof}

\subsubsection{Proof of Theorem \ref{thm_counterex1}}

Suppose that $L$ admits a connected group compactification $(K_1,\tau_1)$ and an index two subgroup $L_2$. For instance, 
we could take 
\[
G = \bZ[1/p] \rtimes \langle p \rangle \qand L = \langle p \rangle \qand N = \bZ[1/p] \qand \Lambda = \bZ,
\]
with $L_2 = \langle p^{2} \rangle$, and $(K_1,\tau_1) = (\bT,\tau_1)$ , where
\[
\tau_1(p^n) = n \log p \mod 1, \quad \textrm{for $n \in \bZ$}.
\] 
Let $K = L/L_2 \times K_1$
and define $\tau_o : L \ra K$ by $\tau_o(l) = (lL_2,\tau_1(l))$. Since $K_1$ is connected, $(K,\tau_o)$ is a group 
compactification of $L$. 

Fix $0 < \eps < 1/2$ and choose a closed $m_{K_1}$-Jordan measurable subset $J_1 \subset K_1$
with $m_{K_1}(J_1) = 2\eps$, equal to the closure of its interior. Pick an element $\delta$ in $L$ such that 
$\delta L_2 \cap L_2 = \emptyset$, and set
\[
I = L_2 \times K_1 \qand J = \delta L_2 \times J_1.
\]
One checks that $I$ and $J$ are closed $m_K$-Jordan measurable sets, equal to the closures of their interiors, and
$m_K(I) = 1/2$ and $m_K(J) = \eps$. Moreover, $J$ has trivial stabilizer in $K$, and $I J = \delta L_2 \times K_1$, which is 
again $m_K$-Jordan measurable, and $m_K(IJ) = 1/2$. Let $A$ and $B$ be as in \eqref{basicform}. Then, by 
Lemma \ref{auxcnt} (i) and \eqref{ineqcntex},
\[
m_K(I) = d^*_G(A) \leq d_G^*(AB) \leq m_K(IJ) = \frac{1}{2} < d^*_G(A) + d^*_G(B) < 1,
\]
which finishes the first part of the conclusion of Theorem \ref{thm_counterex1}. For the second part, note that if $G_o < G$ is any finite-index
subgroup of $G$, then, since $K_1$ is connected, we have either $\overline{\tau(G_o)} = L_2 \times K_1$ or $\overline{\tau(G_o)} = K$. 
Indeed, since $L/L_2$ has two elements, the subgroup $\overline{\tau(G_o)} \cap (L_2 \times K_1)$ must be open in $K_1$ and thus equal 
to $K_1$ since $K_1$ is connected. 

In the case when $\overline{\tau(G_o)} = L_2 \times K_1$, then $\tau(G_o) \cap J = \emptyset$, so by Lemma \ref{auxcnt} (ii,b), 
\[
d^*_G(G_o B) > d_G^*(B) + \frac{1}{[G : G_o]}.
\]
In the case when $\overline{\tau(G_o)} = K$, then $\tau(G_o)J = K$, so by Lemma \ref{auxcnt} (ii,a), we have $G_o B = G$,
which finishes the proof of Theorem \ref{thm_counterex1}.

\subsubsection{Proof of Theorem \ref{thm_counterex2}}
Suppose that $L$ admits a homomorphism $\tau_o : L \ra K \cong \bT$ with dense image; either Example \eqref{cntex1} or 
Example \eqref{cntex2} would do. Pick a closed interval $I \subset \bT$ with $m_{\bT}(I) < 1/3$, such that 
$(I + I) \cap I = \emptyset$; in particular, $0 \notin I$. Let $A$ and $B$ be as in \eqref{basicform} with $I = J$. By Lemma 
\ref{auxcnt} (iii) and (iv), $A$ and $B$ are both spread-out, but $B$ is not contained in a Sturmian set with the same upper 
Banach density as $B$. Furthermore, by Lemma \ref{auxcnt} (i) and \eqref{ineqcntex},
\[
d^*_G(AB) \leq m_K(I+I) = 2 m_K(I) = d^*_G(A) + d^*_G(B) < 1.
\]
Since $A$ is spread-out, Theorem \ref{dens1ubd} tells us that the first inequality cannot be strict, so $(A,B)$ provides the 
counterexample in Theorem \ref{thm_counterex2}.

\appendix
\addtocontents{toc}{\protect\setcounter{tocdepth}{1}}

\section{Invariant means and Furstenberg's Correspondence Principle}
\label{Fur}

\begin{adjustwidth*}{1in}{1in}
\emph{We define amenable groups, invariant means and 
asymptotic densities along F\o lner sequences. We also 
state Furstenberg's well-known Correspondence Principle in 
a slightly unorthodox form, and list some of its useful applications.}
\end{adjustwidth*}

\vspace{0.2cm}

\subsection{Amenable groups and invariant means}
\label{subsection:amenable}

Throughout this section, let $G$ be a countable group. We denote by $\ell^\infty(G)$ the Banach 
space of real-valued bounded functions on $G$, endowed with the uniform norm. If $\lambda$ 
belongs to the dual of $\ell^\infty(G)$, then we shall twice abuse notation when we refer to this 
element. Firstly, we shall identify $\lambda$ with the bounded finitely additive measure $\lambda$ on 
$G$, defined by $\lambda(A) = \lambda(\chi_A)$ for $A \subset G$. Secondly, if $f \in \ell^\infty(G)$,
it will sometimes be convenient to write
\[
\lambda(f) = \int_G f(g) \, d\lambda(g),
\]
although the right hand side is not an integral in the Lebesgue sense. We denote by $\cM_G$ the
weak*-closed and convex set of positive and unital functionals on $\ell^\infty(G)$, so called \emph{means}
on $G$. We say that $\lambda \in \cM_G$ is \emph{invariant} if $\lambda(gA) = \lambda(A)$ for all $g \in G$ and $A \subset G$.
We denote by $\cL_G$ the (possibly empty) set of invariant means on $G$. We say that $G$ is \emph{amenable} 
if $\cL_G$ is non-empty, and we refer to \cite{pa1988} for more information about this class of groups. It suffices for now
to say that every solvable group is amenable, as is every locally finite group and every group of sub-exponential growth.
On the other hand, any group which contains a free subgroup on more than two generators is \emph{not} amenable.

If $G$ is amenable, then $\cL_G$ must contain extreme points by Krein-Milman's Theorem. We denote the set of 
such extreme elements by $\cL_G^{\textrm{ext}}$. The following result, which is quite standard (see for instance \cite{bj2017} for 
a proof), points out an important "ergodicity" property of such means. 

\begin{proposition}[Weak Ergodic Theorem]
\label{wet}
If $\lambda \in \cL_G^{\textrm{ext}}$, then
\[
\int_G \lambda(gA \cap B) \, d\eta(g) = \lambda(A) \, \lambda(B),
\]
for all $A, B \subset G$ and $\eta \in \cL_G$.
\end{proposition}

\subsection{Pointed $G$-spaces}
\label{subsec:pointedG}
Let $G$ be a countable amenable group and suppose that it acts by homeomorphisms on a compact second countable
space $X$. We shall assume that there exists $x_o \in X$ with a dense $G$-orbit, and we refer to the pair 
$(X,x_o)$ as a \emph{pointed $G$-space}. Let $C(X)$ denote the Banach space of real-valued continuous 
functions on $G$ endowed with the uniform norm, and $\cP(X) \subset C(X)^*$ the weak*-closed and convex
set of regular Borel probability measures on $X$. We say that $\mu \in \cP(X)$ is \emph{invariant} if 
$\mu(gB) = \mu(B)$ for all $g \in G$ and every Borel set $B \subset X$, and write $\cP_G(X)$ for the 
set of invariant probability measures. We see that there is a positive and unital linear map 
$S_{x_o} : C(X) \ra \ell^\infty(G)$ defined by $(S_{x_o}f)(g) = f(gx_o)$, which intertwines the left-regular representations
of $G$ on $C(X)$ and $\ell^\infty(G)$.  It readily
follows that $S_{x_o}^*(\cL_G) \subset \cP_G(X)$. In particular, since $G$ is amenable, $\cP_G(X)$ is always 
non-empty. By Krein-Milman's Theorem, the set $\cP^{\textrm{erg}}_G(X)$ of extreme points is then non-empty, and 
it turns out that it coincides with the set of \emph{ergodic} measures in $\cP_G(X)$ (see for instance Theorem 4.4. in \cite{eiwa2015}).

If $\lambda \in \cL_G$ and $\mu = S_{x_o}^*\lambda$, then for every \emph{clopen} set $U \subset X$, the indicator 
function $\chi_U$ is continuous on $X$, and thus we have $\lambda(U_{x_o}) = \lambda(S_{x_o}\chi_U) = \mu(U)$. 
However, this observation is only useful for \emph{disconnected} spaces; as many of the $G$-spaces that we will work 
with are connected, it will be useful to know if this type of identity holds for more general classes of sets (that some regularity on $U$ has
to be assumed is clear already from easy examples). The following lemma provides some useful answers in this 
direction. Recall that if 
$\mu$ is a regular Borel probability measure on $X$, then a set $U \subset X$ is \emph{$\mu$-Jordan measurable} if 
$\mu(\overline{U}) = \mu(U^o)$, where $\overline{U}$ denotes the closure of $U$ and $U^o$ denotes the interior of $U$.

\begin{lemma}
\label{somecons}
Let $\lambda \in \cL_G$ and $\mu = S_{x_o}^*\lambda$. 
\begin{enumerate}
\item If $U \subset X$ is $\mu$-Jordan measurable, then $\lambda(U_{x_o}) = \mu(U)$.
\item If $U, V \subset X$ are $\mu$-Jordan measurable, and $\mu(U \setminus V) = 0$, 
then $\lambda(U_{x_o} \cap V_{x_o}) = \mu(U)$.
\item If $U \subset X$ is open, then $\lambda(U_{x_o}) \geq \mu(U)$.
\end{enumerate}
In particular, if $A \subset G$ is any subset, and $U \subset X$ is clopen, then 
\[
\lambda(AU_{x_o}) \geq \mu(AU) \qand \lambda(U_{x_o}) = \mu(U).
\]
\end{lemma}

\begin{proof}
$(i)$ If $U$ is $\mu$-Jordan measurable, then by Proposition 2.3.3 in \cite{wi2007}, there exist, for every $\eps > 0$, continuous
functions $f_{-}$ and $f_{+}$ on $X$ such that
\[
f_{-} \leq \chi_U \leq f_{+} \qand \mu(f_{+} - f_{-}) \leq \eps,
\] 
and thus 
\[
\mu(f_{-}) = \lambda(S_{x_o}f_{-}) \leq \lambda(U_{x_o}) = 
\lambda(S_{x_o}\chi_U) \leq \lambda(S_{x_o} f_{+}) = \mu(f_{+}),
\]
whence $0 \leq \mu(f_{+}) - \lambda(U_{x_o}) \leq \eps$. By letting $\eps \searrow 0$, we conclude that
$\mu(U) = \lambda(U_{x_o})$.

$(ii)$ One checks that $U \cap V$ is $\mu$-Jordan measurable, and thus by $(i)$, 
\[
\lambda(U_{x_o} \cap V_{x_o}) = \lambda((U \cap V)_{x_o}) = \mu(U \cap V) = \mu(U).
\]
$(iii)$ See Lemma 2.1 in \cite{bj2017_2}.

\end{proof}

In the case when $\cP_G(X) = \{\mu\}$, then $\mu$ is of course extremal in $\cP_G(X)$, and thus ergodic, and for 
every $\mu$-Jordan measurable subset $U \subset X$ and $\lambda \in \cL_G$, we have $\lambda(U_{x_o}) = \mu(U)$, 
no matter if $\lambda$ is extremal in $\cL_G$ or not. In particular, let $(K,L,\tau)$ be an isometric $G$-action. Then, by
Lemma \ref{basic_isomG}, the unique $K$-invariant Borel probability measure $m_{K/L}$ is also the unique $G$-invariant 
probability measure on $X = K/L$. We conclude:

\begin{corollary}
\label{cor_bohrset}
If $I \subset K/L$ is an $m_{K/L}$-Jordan measurable set, then $\lambda(I_t) = m_{K/L}(I)$ for every 
$t \in K/L$ and $\lambda \in \cL_G$. In particular, if $L$ is trivial, then $\lambda(\tau^{-1}(I)) = m_K(I)$
for every $\lambda \in \cL_G$.
\end{corollary}

\subsection{Furstenberg's Correspondence Principle}

We now state Furstenberg's Correspondence Principle in terms of the transpose map $S_{x_o}^*$ above. 
This formulation is perhaps somewhat unorthodox, but can be readily proved along the same lines as in 
Furstenberg's seminal paper \cite{fu1977} where this principle first appeared. A detailed proof of a slightly more 
general statement in the language below can be found in \cite{bj2017_2}.

\begin{proposition}[Furstenberg's Correspondence Principle]
\label{onto}
The map $S_{x_o}^* : \cL_G \ra \cP_G(X)$ 
\begin{enumerate}
\item is affine, weak*-continuous and onto.
\item maps $\cL_G^{\textrm{ext}}$ onto $\cP_G^{\textrm{erg}}(X)$.
\end{enumerate}
\end{proposition}

We stress that it is not at all automatic for weak*-continuous affine maps between weak*-closed and convex sets to map
extreme points to extreme points; indeed, consider the unit square $[0,1]^2$ in $\bR^2$, and the linear map which projects
it onto one of its diagonals. The two corners which are not touched by the diagonal are certainly extreme points of the 
square but will be mapped to midpoint of the diagonal, which is not extreme anymore. 

\subsection{F\o lner sequences and densities}
\label{subsec:folner}
The notions and results in this subsection are well-known, and we only include a brief 
discussion for completeness, and to make referencing easier. 

\begin{definition}
A sequence $(F_n)$ of finite subsets of $G$ is \emph{F\o lner} if
\begin{equation}
\label{defFolner}
\lim_n \frac{|F_n \bigtriangleup g F_n|}{|F_n|} = 0, \quad \textrm{for all $g \in G$}.
\end{equation}
If $(F_n)$ is a F\o lner sequence in $G$, then we define the \emph{upper} and 
\emph{lower asymptotic density} of a subset $A \subset G$ along $(F_n)$ by
\begin{equation}
\label{defasympdens}
\overline{d}_{(F_n)}(A) = \varlimsup_n \, \frac{|A \cap F_n|}{|F_n|}
\qand
\underline{d}_{(F_n)}(A) = \varliminf_n \, \frac{|A \cap F_n|}{|F_n|}
\end{equation}
respectively, and the \emph{upper} and \emph{lower Banach densities} by
\begin{equation}
\label{defBanachdens}
d^*(A) = \sup\big\{ \overline{d}_{(F_n)}(A) \, : \, \textrm{$(F_n)$ F\o lner} \big\}
\qand
d_*(A) = \inf\big\{ \underline{d}_{(F_n)}(A) \, : \, \textrm{$(F_n)$ F\o lner} \big\}
\end{equation}
respectively.
\end{definition}

We note that every F\o lner sequence $(F_n)$ naturally gives rise to invariant means. Indeed, consider the
sequence $(\lambda_n)$ in $\cM_G$ defined by
\[
\lambda_n(f) = \frac{1}{|F_n|} \sum_{g \in F_n} f(g), \quad \textrm{for $f \in \ell^\infty(G)$}.
\]
It readily follows from the F\o lner condition \eqref{defFolner} that any weak*-accumulation point 
of $(\lambda_n)$ is invariant. In particular, for any $A \subset G$, there are $\overline{\lambda}, \underline{\lambda}$ 
such that
\begin{equation}
\label{asympdensmean}
\overline{d}_{(F_n)}(A) = \overline{\lambda}(A)
\qand
\underline{d}_{(F_n)}(A) = \underline{\lambda}(A).
\end{equation}
Of course, for a subset $B \subset G$ different from $A$, we can only be sure of the inequalities 
\begin{equation}
\label{besureof}
\overline{\lambda}(B) \leq \overline{d}_{(F_n)}(B)
\qand
\underline{\lambda}(B) \geq \underline{d}_{(F_n)}(B).
\end{equation}
The following proposition is well-known to experts, but hard to find a good reference for, so we supply
a proof here. 
\begin{proposition}
\label{densextmean}
For every $A \subset G$,
\[
d^*(A) = \sup\big\{ \lambda(A) \, : \, \lambda \in \cL_G \big\}
\qand
d_*(A) = \inf\big\{ \lambda(A) \, : \, \lambda \in \cL_G \big\},
\]
and there are extreme $\lambda_{+}, \lambda_{-}$ in $\cL_G$ such that $d^*(A) = \lambda_{+}(A)$
and $d_*(A) = \lambda_{-}(A)$. 
\end{proposition}

\begin{proof}
Assuming the identities for $d^*$ and $d_*$, the second assertion is immediate from the fact that the map 
$\lambda \mapsto \lambda(A)$ is weak*-continuous and affine on $\cL_G$, and such maps always attain 
their minima and maxima at extreme points. 

Concerning the identities, we first note that \eqref{asympdensmean} implies that 
\[
d^*(A) \leq \sup\big\{ \lambda(A) \, : \, \lambda \in \cL_G \big\}
\qand
d_*(A) \geq \inf\big\{ \lambda(A) \, : \, \lambda \in \cL_G \big\}.
\]
Let us prove that the first inequality is in fact an identity; the second inequality can be treated completely 
analogously. 
Pick an extreme $\lambda$ at which the supremum above is realized, and denote by $(X,x_o)$ the $G$-hull
associated to the set $A$, as in Subsection \ref{subsec:hull}. Abusing notation, we can find a \emph{clopen} 
subset $A \subset X$ such
that our set in $G$ can be represented as $A_{x_o}$. Let $\mu = S_{x_o}^*\lambda$; by Proposition 
\ref{onto}, $\mu$ is ergodic and $\lambda(A_{x_o}) = \mu(A)$. By the strong mean Ergodic Theorem (see e.g. \cite{bjfi2018}), the averages
\[
\lim_n \frac{1}{|F_n|} \sum_{g \in F_n} \chi_A(gx) = \lim_{n} \frac{|A_x \cap F_n|}{|F_n|}, \quad \textrm{for $x \in X$},
\]
converge in $L^2(X,\mu)$ to the constant function $\mu(A)$, whence, upon passing to a sub-sequence $(n_k)$, $\mu$-almost surely
to $\mu(A)$. Pick $x \in X$ for which this sub-sequence converges. Since $A$ is clopen, we can find $(g_{n_k})$ such that
$A_{g_{n_k}.x_o} \cap F_{n_k} = A_x \cap F_{n_k}$ for every $k$, and thus
\[
|A_{x_o} \cap F_{n_k} g_{n_k}| = |A_{g_{n_k}.x_o} \cap F_{n_k}| = |A_x \cap F_{n_k}|,
\]
which shows that
\[
d^*(A_{x_o}) \geq \lim_k \frac{|A_{x_o} \cap F_{n_k} g_{n_k}|}{|F_{n_k}|} = \mu(A),
\]
where the  inequality follows from the fact that $(F_{n_k} g_{n_k})$ is a F\o lner sequence in $G$.
\end{proof}

\subsection{Thickness and syndeticity}

Recall that a subset $A \subset G$ is \emph{thick} if for every finite subset $F \subset G$ there is $g \in G$
such that $Fg \subset A$, and \emph{syndetic} if there exists a finite set $F \subset G$ such that $FA = G$.
For a proof of the following 
well-known density characterizations of thick and syndetic sets, see for instance Subsections 2.5 and 2.6 in \cite{bjfi2015}.

\begin{lemma}
\label{thicksyn}
Let $G$ be a countable amenable group, and let $A \subset G$. Then, 
\begin{enumerate}
\item $A$ is thick $\iff$ $d^*(A) = 1$.
\item $A$ is syndetic $\iff$ $d_*(A) > 0$.
\end{enumerate}
\end{lemma}

\section{Generalities on group compactifications}
\label{sec:gencomp}
\begin{adjustwidth*}{1in}{1in}
\emph{We collect here some basic facts about group compactifications of 
countable groups that will be used in some of our proofs.}
\end{adjustwidth*}

\vspace{0.2cm}

Let $G$ be a countable group. We say that $(K,\tau)$ is a \emph{group compactification} of $G$ if 
$K$ is a compact Hausdorff group and $\tau : G \ra K$ is a homomorphism with dense image. We 
stress that we do \emph{not} assume that $\tau$ is injective. We shall always denote the (unique)
Haar probability measure on $K$ by $m_K$ and the identity element in $K$ by $e_K$.

As it turns out, many group-theoretical properties of $G$ can be transported to topological properties 
of $K$. We record here some instances of this phenomenon.

\begin{lemma}
\label{compactify}
If $(K,\tau)$ is a group compactification of $G$, then
\begin{enumerate}
\item $G$ amenable $\implies$ $K^o$ is abelian.
\item $G$ has no non-trivial finite index subgroups $\implies$ $K$ is connected.
\item $G$ is a finitely generated torsion group $\implies$ $K$ is totally disconnected.
\item $G$ is a finitely generated simple group $\implies$ $K$ is trivial.
\end{enumerate}
\end{lemma}

\begin{proof}
(i) See the Appendix in \cite{bj2017}. 

(ii) If $K$ is not connected, then there is a non-trivial proper open subgroup $U$ of $K$. Since $K$ is compact, $U$ must
have finite index in $K$, and thus $G_o = \tau^{-1}(U)$ has finite index in $G$. 

(iii) By Corollary 2.36 in \cite{homo2013}, 
we can find a net $(N_\alpha)$ of closed normal subgroups of $K$ and integers $(n_\alpha)$ such that
\[
\bigcap_\alpha N_\alpha = \{e_K\} \qand K_\alpha := K/N_\alpha \overset{\iota_\alpha} \hookrightarrow U(n_\alpha),
\]
where $U(n)$ denotes the unitary group in dimension $n$, and $\iota_\alpha$ is injective for every $\alpha$. Note that 
for every $\alpha$, the subgroup $\Gamma_\alpha = \iota_\alpha \circ \tau(G)$ of the linear group $U(n_\alpha)$ is 
finitely generated and torsion. By Jordan-Schur's Theorem, these properties imply that $\Gamma_\alpha$ is finite,
whence $K_\alpha$ is finite, and thus $N_\alpha$ must be open in $K$ for every $\alpha$. Since the intersections of all 
$N_\alpha$ is trivial, $K$ is totally disconnected. 

(iv) First note that (ii) implies that $K$ must be connected, and thus 
Peter-Weyl's Theorem shows that if $K$ is non-trivial, then it admits a non-trivial compact and connected Lie group $K'$ 
as a quotient group. Since $G$ is simple, the composition of $\tau$ with this quotient map is still injective (otherwise the 
kernel would be a non-trivial normal subgroup of $G$). In particular, $G$ can be viewed as a finitely generated 
subgroup of $K'$. However, Malcev's Theorem now says that any such subgroup must be residually finite, and thus far from simple, 
which leads us to conclude that $K$ is trivial.
\end{proof}

The next lemma contains some auxiliary observations about pull-backs of sets in a group compactification $(K,\tau)$ of
a countable group $G$.

\begin{lemma}
\label{stuff}
Let $(K,\tau)$ be a compactification of $G$, and fix a thick subset $T \subset G$, a non-empty
open subset $U \subset K$ and a finite-index subgroup $G_o < G$. Then,
\begin{enumerate}
\item for every $s \in G$, the set $\tau^{-1}(U) \cap G_os$ is non-empty iff it is syndetic. 
\item $G_o(\tau^{-1}(U) \cap T) = \tau^{-1}(\tau(G_o) U)$.
\item if $K$ is connected, then $\overline{\tau(Q \cap T)} = K$ for every non-empty right $G_o$-invariant set $Q \subset G$.
\end{enumerate}
\end{lemma}

\begin{proof}
(i) Note that the closed subgroup $H = \overline{\tau(G_o)} < K$ has finite index, hence open. We note
that $\tau^{-1}(U) \cap G_os = \tau^{-1} (U\tau(s)^{-1} \cap H)s$, and thus, if this set is non-empty, then 
$V = U \tau(s)^{-1}  \cap H$ is a non-empty open subset of $K$. Since $\tau(G)$ is dense in $K$, there is a
finite set $F \subset G$ such that $\tau(F)V = K$, whence $F(\tau^{-1}(U) \cap G_o s) = \tau^{-1}(FV)s = G$,
which shows that $\tau^{-1}(U) \cap G_o s$ is syndetic in $G$. 

(ii) Set $D = \tau^{-1}(U)$, and define 
$D_{+} = \{ s \in G_o \backslash G \, : \, D \cap G_o s \neq \emptyset \big\}$. Note that 
\[
G_o(D \cap T) = \bigsqcup_{s \in D_{+}} G_o((D \cap G_o s) \cap T).
\]
By (i), if $s \in D_{+}$, then $D \cap G_o s$ is in fact syndetic in $G$, and thus intersects the thick set $T$ 
non-trivially, whence $G_o((D \cap G_o s) \cap T) = G_o s$ for all $s \in D_{+}$, which finishes the proof.

(iii) Fix an open identity neighborhood $V$ in $K$ and an exhaustion $(F_n)$ of finite subsets of $G$. Since
$T$ is thick, we can find a sequence $(g_n)$ such that $F_n g_n \subset T$ for all $n$. Since $Q$ is right 
$G_o$-invariant and $K$ compact, we may pass to further sub-sequence (or sub-net, if $K$ is not sequentially compact), 
so that for some $g \in G$ and $t \in K$ we have
\[
Q g_n^{-1} = Q g^{-1} \qand \tau(g_n)^{-1} \in Vt
\]
for all $n$, and thus
\[
\tau(Q \cap T) V  \supset \tau(Q \cap F_n g_n) \tau(g_n)^{-1} t^{-1} = \tau(Q g^{-1} \cap F_n) t^{-1}
\]
for all $n$, whence $\tau(Q \cap T) V \supset \tau(Qg^{-1})t^{-1}$. Since $K$ is connected, $\tau(G_o)$ is dense
in $K$, and thus $\tau(Qg^{-1})$ is dense as well. Since $V$ is arbitrary, we conclude that $\tau(Q \cap T)$ is dense.  
\end{proof}

\section{Peculiar sumsets in $\bZ$ relative to the F\o lner sequence $\{[-n,n]\}$}
\label{sec:pec}

\begin{adjustwidth*}{1in}{1in}
\emph{We show that different attempts to weaken the assumptions in Theorem 
\ref{dens1} and Theorem \ref{dens2} fail, already
for $G = (\bZ,+)$ and the F\o lner sequence $F_n = [-n,n]$.}
\end{adjustwidth*}

\vspace{0.2cm}

To keep things simple, let us in this appendix only focus on the group $G = (\bZ,+)$, the F\o lner sequence $F_n = [-n,n]$ in $G$ 
and its associated \emph{lower} asymptotic density, which we here denote by
\[
\underline{d}(A) = \varliminf_{n \ra \infty} \frac{|A \cap [-n,n]|}{2n+1}, \quad \textrm{for $A \subset \bZ$}.
\]
Our two first examples concern attempts to weaken the hypotheses of Theorem \ref{dens1}, while our third and fourth example deal with failed
conjectural strengthenings of Theorem \ref{dens2}. In each example, the weakened assumption is marked in \textsc{capital} 
letters.

\begin{proposition}
\label{attempt1}
There exist $A, B \subset \bZ$ such that 
\begin{enumerate}
\item $A$ is not contained in a proper periodic set,
\item $B$ is syndetic,
\item $A + B$ is \textsc{thick},
\end{enumerate}
and $\underline{d}(A+B) < d^*(A) + \underline{d}(B) < 1$.
\end{proposition}

\begin{proposition}
\label{attempt2}
There exist $A, B \subset \bZ$  such that 
\begin{enumerate}
\item $A$ is not contained in a proper periodic set,
\item $B$ is \textsc{not syndetic}, but $\underline{d}(B) > 0$,
\item $A + B$ is not thick,
\end{enumerate}
and $\underline{d}(A+B) < d^*(A) + \underline{d}(B) < 1$.
\end{proposition}

\begin{proposition}
\label{attempt3}
There exist $A, B \subset \bZ$ such that
\begin{enumerate}
\item $A$ is spread-out and not contained in a Sturmian set with the same upper Banach density as $A$,
\item $B$ is \textsc{not syndetic}, but $\underline{d}(B) > 0$,
\item $A + B$ does not contain a piecewise periodic set,
\end{enumerate}
and $\underline{d}(A+B) = d^*(A) + \underline{d}(B) < 1$.
\end{proposition}

\begin{proposition}
\label{attempt4}
There exist $A, B \subset \bZ$ such that
\begin{enumerate}
\item $A$ is spread-out and not contained in a Sturmian set with the same upper Banach density as $A$,
\item $B$ is is syndetic,
\item $A + B$ is \textsc{thick},
\end{enumerate}
with $\underline{d}(A+B) = d^*(A) + \underline{d}(B) < 1$.
\end{proposition}

The examples above are constructed by similar procedures, so we will discuss them in parallel. We start
by fixing an irrational $\alpha \in \bT = \bR/\bZ$. Given a proper closed interval $I \subset \bT$, we shall 
write 
\[
C_I = \big\{ n \in \bZ \, : \, n\alpha \in I \big\}.
\]
In Proposition \ref{attempt1}, we pick $I \subset \bT$ with $m_{\bT}(I) < 1/3$, and set $A = C_I \cap \bN$ and 
$B = C_I \cup \bN$. In Proposition \ref{attempt2}, we pick $I \subset \bT$ with $m_{\bT}(I) < 1/2$, and set 
$A = B = C_I \cap \bN$. In Proposition \ref{attempt3}, we pick $I \subset \bT$ such that 
$(I + I) \cap (I+n\alpha) = \emptyset$ for some integer $n$, and set $A = (C_I \cap \bN) \cup \{n\}$ and
 $B = C_I \cap \bN$. The matter of verifying that these choices indeed lead to the examples in the propositions
 is entirely routine, and left to the reader. Proposition \ref{attempt4} is more involved. To construct our example here, 
 we first need to produce a thick set $T \subset \bN$ such that 
 \[
 \underline{d}_{([1,n])}(T) = \frac{1}{10} \qand \underline{d}_{([1,n])}(T+T) = \frac{2}{10},
 \]
 with the property that the sequence $F_n = [1,n] \setminus (T + T)$ is F\o lner. This is tedious, but still a matter of 
 routine. We now choose $I \subset \bT$ with $m_{\bT}(I) = 4/9$ such that for some $m \in \bZ$,
 \[
 m_{\bT}((I+I) \cup (I + m\alpha)) = \big(2 + \frac{1}{24}\big)m_{\bT}(I).
 \]
 The exact numbers here are not so important; the construction has some wiggle room. Once $T$, $I$ and $m$ have
 been produced, we set $A = (C_I \cap T) \cup \{m\}$ and $B = C_I \cup T$, and note that $B$ is syndetic \emph{and}
 thick, so $A+B$ is thick as well.  To check the remaining properties in Proposition \ref{attempt4} is again a matter of
 routine.

\subsection*{Acknowledgments}
The first author has benefited enormously from discussions with Benjy Weiss during the preparation of this manuscript, and
it is a pleasure to thank him for sharing his many insights. 
The authors would also like to acknowledge the great impact that the many conversations with Eli Glasner and John Griesmer have had 
on the work. Furthermore, we have had interesting, enlightening and inspiring discussions with Mathias Beiglb\"ock, Vitaly
Bergelson, Manfred Einsiedler, Hillel Furstenberg, Elon Lindenstrauss, Fedja Nazarov, Amos Nevo, Imre Ruzsa, Omri Sarig
and Jean-Paul Thouvenot. Finally, we would like to thank the referee for a very careful reading of this paper.

The authors started this paper in 2009 at Ohio State University, and continued working on it at University of Wisconsin, Hebrew
University in Jerusalem, Weizmann Institute, IHP Paris, KTH Stockholm, ETH Z\"urich, Chalmers University in Gothenburg and
University of Sydney. Our deepest gratitude goes out to these places for their hospitality.





\begin{dajauthors}
\begin{authorinfo}[michael]
  Michael Bj\"orklund\\
  Department of Mathematics, Chalmers\\ 
  Gothenburg, Sweden\\
  micbjo\imageat{}chalmers\imagedot{}se \\
  \url{http://www.math.chalmers.se/~micbjo/}
\end{authorinfo}
\begin{authorinfo}[sasha]
  Alexander Fish\\
  School of Mathematics and Statistics, University of Sydney\\
  Sydney, Australia\\
  alexander.fish\imageat{}sydney\imagedot{}edu\imagedot{}au \\
  \url{http://www.maths.usyd.edu.au/u/afish/}
\end{authorinfo}
\end{dajauthors}

\end{document}